\documentclass[11pt]{amsart}
\usepackage{amsthm}
\usepackage{amssymb}
\usepackage{amsmath}
\newtheorem{lemma}{Lemma}[section]
\newtheorem{theorem}[lemma]{Theorem}

\newtheorem{proposition}{Proposition}[section]

\newtheorem{corollary}{Corollary}[section]
\newtheorem*{corollary*}{Corollary}

\numberwithin{equation}{section}

\gdef\myletter{}

\let\savetheequation\theequation
\def\theequation{\savetheequation\myletter}

\newcommand{\CC}{{\mathbb C}}

\newcommand{\RR}{{\mathbb R}}

\newcommand{\PP}{{\mathbb P}}

\renewcommand{\Im}{\mbox{Im}}

\renewcommand{\Re}{\mbox{Re}}

\def \bar{\overline}
\def \hat{\widehat}

\begin{document}

\vskip 3mm

\title[Exterior Monge-Amp\`ere ]{\bf Exterior Monge-Amp\`ere Solutions  }  

\author{D. Burns, N. Levenberg and S. Ma'u}{\thanks{Supported in part by NSF grants DMS-0104047 and DMS-0514070 (DB), and by a New Zealand Science \&\ Technology Post-Doctoral fellowship (SM)}}

\address{Univ. of Michigan, Ann Arbor, MI 48109-1109 USA}  
\email{dburns@umich.edu}

\address{Indiana University, Bloomington, IN 47405 USA}

\email{nlevenbe@indiana.edu}

\address{Indiana University, Bloomington, IN 47405 USA}
\email{sinmau@indiana.edu}
\date{July 5, 2006}

\begin{abstract}
We discuss the Siciak-Zaharjuta extremal function of a real convex
body in $\mathbb{C}^n$, a solution of the homogeneous complex
Monge-Amp\`ere equation on the exterior of the convex body. We
determine several conditions under which a foliation by holomorphic
curves can be found in the complement of the convex body along which
the extremal function is harmonic. We study a variational problem for
holomorphic disks in projective space passing through prescribed
points at infinity.  The extremal curves are all complex quadric
curves, and the geometry of such curves allows for the determination
of the leaves of the foliation by simple geometric criteria. As a
byproduct we encounter a new invariant of an exterior domain, the Robin
indicatrix, which is in some cases the dual of the Kobayashi
indicatrix for a bounded domain.
\end{abstract}

\maketitle

\section{\bf Introduction.}
	\label{sec:intro}

\vskip 3mm

For a function $u$ of class $C^2$ on a domain $\Omega\subset \CC^n$, the complex Monge-Amp\`ere operator applied to $u$ is
$$(dd^cu)^n:= i\partial \bar \partial u \wedge \cdots \wedge i\partial \bar \partial u$$
$$ =4^n n! \det [\frac{\partial^2 u}{\partial z_j \partial \bar z_k}]_{j,k=1,...,n}\cdot (i/2)^n\prod_{j=1}^n dz_j\wedge d\bar z_j.$$
Thus if $u$ is plurisubharmonic (psh) and satisfies the homogeneous complex Monge-Amp\`ere equation $(dd^cu)^n=0$ in $\Omega$, then at each point of $\Omega$ the complex Hessian of $u$ has a zero eigenvalue. For a psh function $u$ which is only locally bounded, $(dd^cu)^n$ is well-defined as a positive measure (cf., \cite{bt}). 

Given a bounded set $E\subset \CC^n$, the 
Siciak-Zaharjuta extremal function is defined as
$$V_E(z):=\sup \{u(z):u\in L(\CC^n), \ u\leq 0 \ \hbox{on} \ E\}$$
where $L(\CC^n)$ denotes the class of psh functions $u$ on $\CC^n$ with $u(z)\leq \log ^+{|z|}+ c(u)$, for some real constant $c(u)$. If 
$E$ is compact, this coincides with 
$$\sup \{\frac{1}{{\rm deg} \, p}\log |p(z)|: {\rm deg} \, p>0, \ ||p||_E:=\sup_{z\in E}|p(z)|\leq 1\}.$$
If $E$ is non-pluripolar, the function $V_E^*(z):=\limsup_{\zeta\to z}V_E(\zeta)$ is a locally bounded psh function which satisfies $(dd^cV_E^*)^N=0$ outside of $\bar E$. The 
relevant situations in this paper occur for $E=\overline D$, $D$ bounded, smoothly bounded and strictly lineally convex in $\CC^n$, and
for $E=K\subset \RR^n$ a convex body, that is, a compact, convex set with non-empty interior. In the former case, $\CC^n \setminus \bar D$ is foliated by
holomorphic curves $f(\CC\setminus \bar \triangle)$ ($\triangle$ denoting the unit disk) on which the Monge-Amp\`ere solution $V_{\bar {D}}$ is harmonic
\cite{momm}, and $V_{\bar{D}}$ is smooth on all of $\mathbb{C}^n \setminus D.$ In the latter situation the function $V_K$ is continuous but it is not smooth. For a symmetric convex body ($K=-K$) it was  shown (\cite{lun}, \cite{bar}) 
that there exists a continuous foliation of $\CC^n\setminus K$ by analytic disks on which $V_K$ is harmonic. 

If $K$ is not necessarily symmetric, 
the existence, through each point in $\CC^n \setminus K$, of an analytic disk on which 
$V_K$ is harmonic was demonstrated in \cite{blm}. This was acheived by  
using a decreasing sequence of strictly convex $\sigma$-invariant ($\sigma (z):=\bar z$) open sets $D_j$ such that $\cap D_j = K$, and a normal families argument on
a sequence of foliation curves $f_j(\CC\setminus \triangle)$ for $\CC^n \setminus D_j$.

In this paper we return to the matters discussed in \cite{blm} to examine the 
properties of the Monge-Amp\`ere solution more closely. In particular, our goal is to clarify the situation for exterior Monge-Amp\`ere solutions $V_K$ associated to general real convex bodies. We use a variational description of the limit disks introduced in \cite{dbx} to base the study of the curves on their behavior at infinity in projective space $\CC\PP^n$ containing $\CC^n$. We give criteria for these curves to give a foliation of $\CC^n \setminus K$ which include the case when the variational problem has a unique extremal curve. When these extremal curves are not unique, we can sometimes retrieve a continuous foliation of $\CC\PP^n \setminus K$ by extremals, and even make a canonical choice of this foliation.

We begin in sections 2 and 3 by verifying that
for a bounded, smoothly bounded and strictly lineally convex domain $D$ containing the origin, and its dual, $D'$, the Kelvin transform ${\mathcal K}: D' \setminus \{0\}
\rightarrow \CC^n \setminus \bar{D}$ described by Lempert in \cite{lem1} and utilized in \cite{blm} extends to a diffeomorphism $\hat{{\mathcal K}}: \hat{D}' \rightarrow
\CC\PP^n \setminus \bar{D}$, where $\hat{D}'$ is the blow-up of $D'$ at the origin. This allows us to parametrize the leaves of the foliation of $\CC^n \setminus
\bar D$ by the hyperplane at infinity, $H_{\infty}$, and to set up a variational problem in section 4 for extremal disks $f_c, \ [c]\in H_{\infty}$, in the spirit
of \cite{lem1}, and dual in a sense to the Kobayashi-Royden functional for the infinitesimal Kobayashi metric. In section 3, a key role is played by the second variation of the Kobayashi-Royden functional, already examined in \cite{dbx}. Using this technique, we note in passing that the Kobayashi indicatrix $I_0(D') \subset T'_0(D')$ of a strictly lineally convex $D' \subset \subset \CC^n$ is itself strictly convex (Corollary \ref{cor:KI}). This indicatrix is, in some sense, the best circled approximation of $D'$ (see \cite{lem3}). In direct analogy to this, in section \ref{sec:varprob} we note that the individual extremal disks in our exterior problem may be packaged into an exponential map from an indicatrix, the Robin indicatrix, defined in the holomorphic normal bundle to $H_{\infty}$, to the complement $\CC\PP^n\setminus{\bar{D}}$. In sections 5 and 6 we pass to the limit $\cap D_j = K$, as in \cite{blm}, and use the special form of the analytic disks in this case -- they are complex ellipses -- to give a simple but useful geometric interpretation to the extremals of the variational problem. Section 7 concerns a geometric criterion on $K$ for uniqueness of extremal curves passing through a given $[c] \in H_{\infty}$. In section 8, we verify that if extremal curves $f_c$ are unique for a real convex
body $K$, then we obtain a continuous foliation as in the symmetric case. This occurs, for example, if $K$ is strictly convex, but the general condition is somewhat weaker than this. This foliation property  depends on the geometry of the real ellipses in $K$ given by the intersection of the complex extremal curves with $\RR^n$. In section 9 we show that in
$\RR^2$ for an arbitrary convex body we may get many foliations of the complement of $K$ by constructing appropriate approximations of $K$ by convex bodies $K_j$
with unique extremals. Then, in section 10, we construct one such foliation, albeit by different means, which is, in a natural way, a ``canonical'' foliation. In fact, the set of extremals through a given $[c] \in H_{\infty}$ is parametrized by a convex closed subset of $K$, and the canonical curve $f_c$ is the one parametrized by the barycenter of this set. In the case where $K = -K$, this foliation reduces to that of \cite{lun} and \cite{bar}.   Finally, in section \ref{sec:ris} we relate the Robin indicatrix of an arbitrary convex body in $\RR^n$ to that of a naturally associated symmetric one. 

The present work is limited to domains $K \subset \RR^2$ in sections 9 and 10 due to a topological argument (continuity of intersection) we use to show that our family of extremal disks is actually a foliation (cf., the proof of Proposition 9.1). It is not clear that this is an essential obstruction, and this should be clarified. We hope to discuss the regularity of the foliation described here and the corresponding regularity of the Monge-Amp\`ere solution $V_K$ in a future paper. The regularity depends on a more detailed study of the behavior of the real ellipses in $K$ used in this paper. Moreover, we hope to define certain invariants associated to a convex body in $\RR^n$ and its Robin indicatrix {\em \`a la} Lempert \cite{lem3}.

\vskip 5mm

\section{\bf The Behavior of the Kelvin Transform at the origin in $D'$.}
	\label{sec:KT}

\vskip 3mm

We are given a domain $D$, bounded, smoothly bounded and strictly lineally convex in $\CC^n$, containing the origin. This means for each $a\in \partial D$, the
complex tangent hyperplane
$$T_a^{\CC}(\partial D):=\{\zeta=(\zeta_1,...,\zeta_n)\in {\CC}^n:\sum_{j=1}^n(\zeta_j-a_j)\frac{\partial \rho}{\partial z_j}(a)=0\},$$
where $\rho\in C^2(\overline D)$ is a defining function for $D$, satisfies $T_a^{\CC}(\partial D)\cap \overline D =\{a\}$, and $\rho$ restricted to $T_a^{\CC}(\partial D)$ has a non-degenerate minimum at $a$. We want to consider the
Monge-Amp\`ere solution 
$V_{\bar {D}}$ on the exterior of $\bar{D} \subset \subset \CC^n$. According to Lempert \cite{lem2}, there is a  domain 
$$D':=\{z=(z_1,...,z_n)\in \CC^n: \sum_{j=1}^n z_j p_j \not= 1 \ \hbox{for all} \ p=(p_1,...,p_n)\in D\}$$ containing the origin dual to $D$, so that $D'$ is
smoothly bounded and strictly lineally convex in $\CC^n$. By \cite{lem1}, $D'$ admits a plurisubharmonic Green's function 
\begin{eqnarray*} G(z) &=&G_{D'}(z)\\&:=&\sup\{u(z): u \ \hbox{psh in} \ D', \ u\leq 0, \ u(z)-\log {|z|}=0(1), \ z\to 0\}\end{eqnarray*} for $z \in D'$ which is smooth in $D' \setminus
\{0\}$ and has a logarithmic singularity at $z = 0$. The Kelvin transform is a diffeomorphism $${\mathcal K}: D' \setminus \{0\} \rightarrow \CC^n
\setminus
\bar{D},$$ given by the formula 
\begin{equation}
\label{Kelvin}
{\mathcal K}(z) = \frac{(G_{z_1}(z),\ldots, G_{z_n}(z))}{z_1 G_{z_1}(z) + \ldots + z_n G_{z_n}(z)},
\end{equation}
where $G_{z_j}(z) =
\frac{\partial G}{\partial z_j}(z), \ j = 1,\ldots, n.$ Lempert's formula for $V_{\bar {D}}$ is 
\begin{equation} 
\label{eqn:greensiciak}
V_{\bar {D}} ({\mathcal K}(z)) = - G(z).
\end{equation}
Reference \cite{lem1} provides detailed information
about the nature of the singularity of $G(z)$ at $z = 0$, and in order to be able to use this to understand $V_{\bar {D}}$ better at infinity, we will have to examine the
behavior of the Kelvin transform near $z = 0$. It is useful to replace $D'$ by $\hat{D}'$, the blow-up of $D'$ at the origin. Let $\pi$ be the projection $\pi:
\hat{D}'
\rightarrow D'$, and ${\mathcal E} = \pi^{-1}(0)$ the exceptional divisor. Similarly it is useful to add the divisor at infinity $H_{\infty} = \CC\PP^n \setminus \CC^n$ to
the domain $\CC^n
\setminus \bar{D}$, replacing it by $\CC\PP^n \setminus \bar{D}$.  Let $\hat{{\mathcal K}}: \hat{D}' \setminus {\mathcal E} \rightarrow \CC^n \setminus \bar{D} \subset \CC\PP^n
\setminus \bar{D}$ be the diffeomorphism given by $\hat{{\mathcal K}} = {\mathcal K} \circ \pi$. The exact statement we want is the following.

\begin{theorem}
	\label{th:kelinf}
The map $\hat{{\mathcal K}}$ extends to a smooth diffeomorphism, also denoted $\hat{{\mathcal K}}$,
\begin{displaymath}
\hat{{\mathcal K}}: \hat{D}' \rightarrow \CC\PP^n \setminus \bar{D}.
\end{displaymath}
\end{theorem}

\begin{proof} 
We first show that $\hat{{\mathcal K}}$ extends smoothly across ${\mathcal E}$, sending ${\mathcal E}$ to $H_{\infty}$. That the resulting map is a diffeomorphism is more delicate. 
The problem will be set up in this section, but completed only in the next.

\vskip 2mm

We first have to recall some standard facts about coordinates on $\hat{D}'$ and on $\CC\PP^n$. First of all, note that the exceptional divisor ${\mathcal E}$ is 
identified with $\PP(T_0'D')$, the projective space of holomorphic tangent directions at $0 \in D'$. We can assume without loss of generality that we are trying to
extend $\hat{{\mathcal K}}$ smoothly across the point in ${\mathcal E}$ corresponding to the direction $\frac{\partial}{\partial z_1} \in T'_0D'$. Then there exist canonical coordinates
$(z_1, \eta_2,\ldots, \eta_n)$ locally on $\hat{D}'$, so that $[\frac{\partial}{\partial z_1}] = 0$, and the mapping $\pi$ is given by 

\vskip 2mm

\centerline{$\begin{array}{rcl}
	\pi^{*} z_1 & = & z_1 \\
	\pi^{*} z_j & = & z_1 \cdot \eta_j, \; \;  j = 2,\ldots, n.
\end{array}$}

\vskip 2mm

It follows from \cite{lem1} that 
\begin{equation} 
\label{eqn:pistar}
(\pi^{*}G) (z_1, \eta)= \log |z_1| + H(z_1, \eta),
\end{equation}
where $\eta = (\eta_2,\ldots, \eta_n),$ and $H(z_1, \eta)$ extends smoothly across ${\mathcal E}$.
To calculate $\hat{{\mathcal K}}$ locally at $0$, we note that

\vskip 2mm

\centerline{$\begin{array}{rcl}
	\pi^{*} dz_1 & = & dz_1 \\
	\pi^{*} dz_j & = & \eta_j dz_1 + z_1 d\eta_j, \; \;  j = 2,\ldots, n,
\end{array}$}

\vskip 2mm

\noindent from which we conclude immediately

\vskip 2mm

\centerline{$\begin{array}{rcl}
	 \frac{\partial}{\partial z_1} & = & d\pi_{*}(\frac{\partial}{\partial z_1}  - \frac{1}{z_1} \cdot \sum_{j = 2}^{n} \eta_j \frac{\partial }{\partial \eta_j}) \\
	 	&	&	\\
	 \frac{\partial }{\partial z_j} & = & d\pi_{*} (\frac{1}{z_1} \frac{\partial }{\partial \eta_j} ), \; \;  j = 2,\ldots, n.
\end{array}$}

\vskip 2mm

In turn, we conclude that

$$\pi^{*}\frac{\partial G}{\partial z_1} = \frac{1}{2z_1} + \frac{\partial H}{\partial z_1} - \frac{1}{z_1} (\eta_2 \frac{\partial H}{\partial \eta_2} + 
\ldots + \eta_n \frac{\partial H}{\partial \eta_n}),$$ and $$\pi^{*} \frac{\partial G}{\partial z_j} = \frac{1}{z_1} \frac{\partial H}{\partial \eta_j}, \ j =
2,\ldots, n.$$ As to the denominator of $\hat{{\mathcal K}}$, we obtain $$\pi^{*}(z_1 \frac{\partial G}{\partial z_1} + \ldots + z_n \frac{\partial G}{\partial z_n}) = 1/2 +
z_1 \frac{\partial H}{\partial z_1}.$$ Gathering terms, we get $$\hat{{\mathcal K}}^{*}  z_1 = \frac{1}{z_1} (\frac{1/2 + z_1 \frac{\partial H}{\partial z_1} - \sum_{j =
2}^n
\eta_j \frac{\partial H}{\partial \eta_j}}{1/2 + z_1 \frac{\partial H}{\partial z_1}}),$$ and $$\hat{{\mathcal K}}^{*} z_j = \frac{1}{z_1} (\frac{\frac{\partial
H}{\partial
\eta_j}}{1/2 + z_1 \frac{\partial H}{\partial z_1}}), \ j = 2,\ldots, n.$$ 

To examine the limit behavior of $\hat{{\mathcal K}}$ as $z_1 \rightarrow 0$, i.e., near ${\mathcal E}$, we need affine coordinates in a neighborhood of the hyperplane $H_{\infty}$. 
It will suffice to look at $$w_1 = \frac{1}{z_1},$$ and $$w_j = \frac{z_j}{z_1}, \ j = 2,\ldots, n.$$ In these coordinates, $H_{\infty}$ is given by $w_1 = 0$.
Substituting these into the formulas above gives us $$\hat{{\mathcal K}}^{*}w_1 = z_1 \frac{1/2 + z_1 \frac{\partial H}{\partial z_1}}{1/2 + z_1 \frac{\partial H}{\partial z_1} -
\sum_{j = 2}^n \eta_j \frac{\partial H}{\partial \eta_j}},$$ and $$\hat{{\mathcal K}}^{*}w_j = \frac{\frac{\partial H}{\partial \eta_j}}{1/2 + z_1 \frac{\partial H}{\partial
z_1} - \sum_{j = 2}^n \eta_j \frac{\partial H}{\partial \eta_j}}, \; j = 2,\ldots, n.$$ From this it is clear that $\hat{{\mathcal K}}$ extends smoothly across ${\mathcal E}$ near $z_1
= 0, \ \eta = 0$, sending ${\mathcal E}$ to $H_{\infty}$ and sending specifically $(z_1, \eta) = (0, 0)$ to the point $w = (0, 2\frac{\partial H}{\partial
\eta_2},\ldots,
 2\frac{\partial H}{\partial \eta_n})$ in our coordinates $(w_1,\ldots, w_n)$. More precisely, we will see in the next section that $(z_1, \eta) = (0, 0)$ gets sent to the point $(0,\ldots,0)$. 

\vskip 2mm

To complete the proof of Theorem \ref{th:kelinf}, it will suffice to show that the differential $d\hat{{\mathcal K}}_{*}$ is non-singular at $(z_1, \eta) = (0, 0)$. 
This is because $H_{\infty}$ is compact and simply-connected, and so $\hat{{\mathcal K}}$ would induce a covering and therefore a diffeomorphism of ${\mathcal E}$ to $H_{\infty}$, and
therefore of a neighborhood of ${\mathcal E}$ to a neighborhood of $H_{\infty}$. The diffeomorphism property away from ${\mathcal E}$ was already shown in \cite{lem2}.

\vskip 2mm

Let us now calculate the (real) differential $d\hat{{\mathcal K}}_{*}$ at $(z_1, \eta) = (0, 0)$. 

$$d\hat{{\mathcal K}}_{*}(0, 0) = \left(\begin{array}{cccc|cccc}
\frac{\partial w_1}{\partial z_1} & \frac{\partial w_1}{\partial \eta_2} & \ldots & \frac{\partial w_1}{\partial \eta_n} & \frac{\partial w_1}{\partial \bar{z}_1} & \frac{\partial w_1}{\partial \bar{\eta}_2} & \ldots & \frac{\partial w_1}{\partial \bar{\eta}_n} \\
\cdot & \cdot &\cdot & \cdot & \cdot & \cdot & \cdot & \cdot \\
\frac{\partial w_n}{\partial z_1} & \frac{\partial w_n}{\partial \eta_2} & \ldots & \frac{\partial w_n}{\partial \eta_n} & \frac{\partial w_n}{\partial \bar{z}_1} & \frac{\partial w_n}{\partial \bar{\eta}_2} & \ldots & \frac{\partial w_n}{\partial \bar{\eta}_n} \\ 
	&	&	&	&	&	&	&	\\ \hline
	&	&	&	&	&	&	&	\\	
\frac{\partial \bar{w}_1}{\partial z_1} & \frac{\partial \bar{w}_1}{\partial \eta_2} & \ldots & \frac{\partial \bar{w}_1}{\partial \eta_n} & \frac{\partial \bar{w}_1}{\partial \bar{z}_1} & \frac{\partial \bar{w}_1}{\partial \bar{\eta}_2} & \ldots & \frac{\partial \bar{w}_1}{\partial \bar{\eta}_n} \\
\cdot & \cdot & \cdot & \cdot & \cdot & \cdot & \cdot & \cdot \\
\frac{\partial \bar{w}_n}{\partial z_1} & \frac{\partial \bar{w}_n}{\partial \eta_2} & \ldots & \frac{\partial \bar{w}_n}{\partial \eta_n} & \frac{\partial \bar{w}_n}{\partial \bar{z}_1} & \frac{\partial \bar{w}_n}{\partial \bar{\eta}_2} & \ldots & \frac{\partial \bar{w}_n}{\partial \bar{\eta}_n}
\end{array} \right).$$

At $(0, 0)$, $$\frac{\partial w_1}{\partial z_1} = 1, \ \frac{\partial w_1}{\partial \bar{z}_1} = \frac{\partial w_1}{\partial \eta_j} = \frac{\partial
w_1}{\partial \bar{\eta}_j} = 0, \ j = 2,\ldots, n.$$ Using elementary row and column operations, it suffices, therefore, to show that the matrix
$$A = \left(\begin{array}{c|c}
 \frac{\partial w_i}{\partial \eta_j}  & \frac{\partial w_i}{\partial \bar{\eta}_j} \\

		&		\\ \hline
		&		\\	
 \frac{\partial \bar{w}_i}{\partial \eta_j} & \frac{\partial \bar{w}_i}{\partial \bar{\eta}_j} 
 \end{array} \right), \; i, j = 2,\ldots, n,$$ is non-singular at $(0, 0)$.
Since
$$\frac{\partial w_i}{\partial \eta_j}(0, 0) = 2\bigl [\frac{\partial^2 H}{\partial \eta_i \partial \eta_j} + 2\frac{\partial H}{\partial \eta_i} \frac{\partial H}{\partial \eta_j}\bigr ],$$ and 
$$\frac{\partial w_i}{\partial \bar{\eta}_j}(0, 0) = 2\frac{\partial^2 H}{\partial \eta_i \partial \bar{\eta}_j},$$ and the function $H$ is real-valued, we obtain
$$A = \left(\begin{array}{c|c}
 2\bigl [ \frac{\partial^2 H}{\partial \eta_i \partial \eta_j} + 2\frac{\partial H}{\partial \eta_i} \frac{\partial H}{\partial \eta_j}\bigr ]& 2\frac{\partial^2 H}{\partial \eta_i \partial \bar{\eta}_j}  \\ 
 	&	\\ \hline
	&	\\
 2\frac{\partial^2 H}{\partial \eta_j \partial \bar{\eta}_i} & 2\bigl [\frac{\partial^2 H}{\partial \bar{\eta}_i \partial \bar{\eta}_j} +2 \frac{\partial H}{\partial \bar{\eta}_i} \frac{\partial H}{\partial \bar{\eta}_j}\bigr ]
\end{array} \right), \; i,j = 2,\ldots, n.$$ In the next section we will give a geometric interpretation of the matrix $A$, and show that for the $D'$ under 
consideration here it is non-singular.

\end{proof}

\vskip 5mm

\section{\bf Stability and Convexity of the Kobayashi Indicatrix.}
	\label{sec:stability}
	
\vskip 3mm

In this section we will recall the interpretation of the function $H(0, \eta)$, and show the non-singularity of the matrix $A$ is equivalent to the lineal 
convexity of the Kobayashi indicatrix at $z = 0 \in D'$, which is implied by the stability of the Kobayashi-Royden functional at $0$.

First of all, let $f_{\eta} (\zeta):= f(\zeta; \eta)$ be the Kobayashi-Royden extremal disk $f_{\eta}: \triangle \rightarrow D'$, with $f_{\eta}(0) = 0 \in D'$ and 
with 
$$f'_{\eta}(0) = \lambda_{(1,\eta)}\cdot (1, \eta)=\lambda(f_\eta)\cdot (1, \eta)$$ where $\lambda_{(1,\eta)}=\lambda(f_\eta)$ is maximal over 
\begin{equation}
\label{lambda}
\{\lambda > 0 | \, \mbox{there is} \, f\, \mbox{holomorphic}, \, f: \triangle
\rightarrow D', \ f(0) = 0, \ f'(0) = \lambda (1, \eta)\}.
\end{equation}
Let $\hat{f}_{\eta}$ be the proper transform of $f_{\eta}$, i.e., $\pi \circ  \hat{f}_{\eta}=f_{\eta}$, so that $\hat{f}_{\eta}(0) = (0, \eta),$
and
$\hat{f}'_{\eta}(0) = \lambda_{(1,\eta)} \frac{\partial }{\partial z_1} + v_{\eta}$, where $v_{\eta}$ is tangent to ${\mathcal E}$ at $(0, \eta)$. It follows from \cite{lem1} that
\begin{equation}
\label{eqn:grkob}
(f^{*}_{\eta} G)(\zeta) = \log |\zeta| \ \hbox{for all} \ \eta, 
\end{equation}
and that $f(\zeta; \eta)$ is smooth in $(\zeta, \eta)$. Lifting equation (\ref{eqn:grkob}) to $\hat D'$, we get, using (\ref{eqn:pistar}), 
$$\log |\zeta| = \hat{f}_{\eta}^{*} \circ \pi^{*} G = \log |\hat{f}_{\eta, 1}| + \hat{f}_{\eta}^{*}H,$$ and evaluating at $\zeta = 0$ gives 
\begin{equation}
\label{ho}
H(0, \eta) = -\log
\lambda_{(1,\eta)}.
\end{equation}
But then, by (\ref{lambda}) and the definition of the infinitesimal Kobayashi metric $K_0$ on $T'_0(D')$, we have $K_0(1, \eta)=\lambda_{(1,\eta)}^{-1}$ so that 
$$e^{H(0, \eta)} = K_0(1, \eta).$$ We will
show that $\eta = 0$ is a critical point for $K_0$ restricted to this complex hyperplane through $(1, 0)$, and that its real Hessian is positive definite there.

\begin{theorem}
	\label{th:2.1}
\noindent a) For $j = 2,\ldots, n, \ \frac{\partial H}{\partial \eta_j}(0, 0) = 0.$

\noindent b) For $a = (a_2,\ldots,a_n) \in \CC^{n-1}, \ a \neq 0$, we have $$\frac{d^2}{dt^2} H(0, ta) |_{t=0} = 2\Re\{\sum_{2 \leq i, j \leq n} \frac{\partial^2
H}{\partial \eta_i \partial \eta_j}(0, 0)a_i a_j + \frac{\partial^2 H}{\partial \eta_i \partial \bar{\eta}_j}(0, 0) a_i \bar{a}_j\} > 0.$$

\end{theorem}

\begin{corollary}
	\label{cor:2.1}
	The matrix A above is non-singular.
\end{corollary}

As a result of Theorem \ref{th:2.1}, we see that the hyperplane $\{(1, \eta) | \eta \in \CC^{n-1}\}$ is tangent to the level set 
$\{v: K_0(v) = K_0(1, 0) = e^{H(0, 0)}\},$ and that $K_0$ restricted to that hyperplane is strictly convex near the point $(1, 0)$. Recall that the
Kobayashi indicatrix $I_0(D') := \{v \in T'_0(D') | K_0(v) < 1\}$. Thus we see that $I_0(D')$ is strictly lineally convex. Moreover, since $K_0(tv) = |t| K_0(v)$, we have the following corollary (cf., \cite{lem3}, Remark 16.2).

\begin{corollary}
	\label{cor:KI}
The Kobayashi indicatrix $I_0(D')$ is strictly convex.
\end{corollary}

We turn to the proof of Theorem \ref{th:2.1}.
\begin{proof} 
	\label {pf:2.1}
	We write $\mathcal{O}(\triangle)=\mathcal{O}(\triangle,D')$ for the holomorphic maps $f:\triangle\to D'$. Define first, for $k$ an integer $> 4, 0 < \alpha < 1$, $\mathcal{H}_{(0, (1, \eta))}^{k, \alpha}(D') \subset \mathcal{C}^{k, \alpha}(\bar{\triangle}) 
\cap \mathcal{O}(\triangle)$ by 
$$\mathcal{H}_{(0, (1, \eta))}^{k, \alpha}(D') = \{f | f(\partial \triangle) \subset \partial D', \ f(0) = 0, \ f'(0) = \lambda (1,
\eta), \ \lambda = \lambda(f) > 0\}.$$ 
Following \cite{lem1}, in a neighborhood of $f_{0}(\bar \triangle)$ we may find holomorphic coordinates $w = (w_1,\ldots,
w_n)$ and a defining function $r$ for $D'$ such that $f_{0}(\zeta) =f(\zeta;{0})=(\zeta, 0,\ldots, 0), \ \zeta \in \triangle$, and for $\zeta \in \partial \triangle$ we have:

\vskip 3mm

\begin{equation}
\label{eqn:rnormalizations}
\begin{array}{rl}
i) & r_{w_1}(f_{0}(\zeta)) = \bar{\zeta}, \ r_{w_j}(f_{0}(\zeta)) =  0, \  j=2,...,n  \\
	&	\\
ii) & r_{w_1, w_j}(f_{0}(\zeta)) = 0,  \ j = 1,\ldots, n \\
	&	\\
iii) & r_{w_i, \bar{w}_j}(f_{0}(\zeta)) = \delta_{i, \bar{j}}  \\
	&	\\
iv) & \mbox{for} \, a = (a_2,\ldots, a_n) \in \CC^{n-1}, \ a \neq 0,\\
	&	\\
	& \sum_{i, j = 2}^{n} \Re\{r_{w_i, \bar{w}_j}(f_{0}(\zeta)) a_i \bar{a}_j + r_{w_i, w_j}(f_{0}(\zeta)) a_i a_j\} > 0.
\end{array}
\end{equation}

\vskip 3mm

We can now write, for $f=(f_1,...,f_n)\in \mathcal{H}_{(0, (1, \eta))}^{k, \alpha}(D')$, $$\lambda = \lambda(f) = \frac{1}{2\pi i} \int_{\partial \triangle}
\frac{f_1(\zeta)}{\zeta^2} d\zeta  = \frac{1}{2\pi} \int_0^{2\pi} \frac{f_1(\zeta)}{\zeta} d\theta.$$ Since $\lambda$ is real, we can write this $$\lambda(f) =
\Re\{\frac{1}{2\pi} \int_0^{2\pi}
\frac{f_1(\zeta)}{\zeta} d\theta\} = \frac{1}{2\pi} \int_0^{2\pi} \Re\{\frac{f_1(\zeta)}{\zeta}\} d\theta.$$ Let $f(\zeta; \eta)$ be the Kobayashi-Royden
extremal disk for $\eta$, and let $$\dot{f} =
\dot{f}(\zeta; a) = \frac{d}{dt} f(\zeta; ta) |_{t=0}$$ be the corresponding variation vector.  As in \cite{pol}, we differentiate at $t = 0$ to get
$$\frac{d}{dt} \lambda(f(\zeta; ta)) |_{t = 0} = \frac{1}{2\pi} \int_0^{2\pi} \Re\{\frac{\dot{f}_1(\zeta; a)}{\zeta}\} d\theta.$$ To determine $\dot{f}(\zeta; a)$, 
we differentiate $r(f(\zeta; ta)) \equiv 0$, at $t = 0$ and all $\zeta \in \partial \triangle$ and use $i)$:
$$0 = \Re\{\sum_{j = 1}^n r_{w_j}(f_{0}(\zeta))\cdot \dot{f}_j(\zeta; a)\} = \Re \; \bar{\zeta} \, \dot{f}_1(\zeta; a) = \Re \; \frac{\dot{f}_1(\zeta; a)}{\zeta},$$
from which we conclude $\frac{d}{dt} \lambda(1,ta) |_{t = 0} = 0$. This says that
$$\sum_{j=2}^n\frac{\partial \lambda}{\partial \eta_j}(0)a_j=0$$
for all $a$; taking successively $a=(1,...,0), \ a=(0,1,...,0),...$, and $a=(0,...,1)$, we conclude that $\frac{\partial \lambda}{\partial \eta_j}(0)=0$ for
$j=2,...,n$. Using $H(0, \eta) = -\log \lambda_{(1,\eta)}$, this proves Theorem \ref{th:2.1}, part a).

As to part b) of the theorem, we have
$$\frac{d^2}{dt^2} \lambda_{(1,ta)} |_{t = 0} = \frac{1}{2\pi} \int_0^{2\pi} \Re\{\frac{\ddot{f}_1(\zeta; a)}{\zeta}\} d\theta.$$ Differentiating 
$r(f(\zeta; ta)) \equiv 0$ twice at $t = 0$ and using $ii)$ and $iii)$, we obtain 
$$0 = \Re\{\bar{\zeta} \ddot{f}_1(\zeta; a)
+ 2\sum_{i, j = 2}^n \bigl[ r_{w_i, \bar{w}_j} (f_{0}(\zeta))\cdot [\dot{f}_i (\zeta; a) \bar{\dot{f}_j(\zeta; a)}]$$ 
$$+ r_{w_i,
w_j} (f_{0}(\zeta))\cdot [\dot{f}_i(\zeta; a) \dot{f}_j(\zeta; a)]\bigr] \}+|\dot{f}_1(\zeta; a)|^2,$$ and we conclude 
$$\Re \{\frac{\ddot{f}_1(\zeta; a)}{\zeta} \}$$ 
$$= -2\Re\{\sum_{i, j = 2}^n r_{w_i, \bar{w}_j} (f_{0}(\zeta))\cdot [\dot{f}_i (\zeta; a) \bar{\dot{f}_j(\zeta; a)}] + 
r_{w_i,
w_j} (f_{0}(\zeta))\cdot [\dot{f}_i(\zeta; a) \dot{f}_j(\zeta; a)]\}$$
$$-|\dot{f}_1(\zeta; a)|^2,$$ and so from $iv)$
\begin{equation}
\label{eqn:2varsign}
\frac{d^2}{dt^2} \lambda_{(1,ta)} |_{t = 0}= \frac{-1}{\pi} \int_0^{2\pi} \Re\{\sum_{i, j = 2}^n r_{w_i, \bar{w}_j} (f_{0}(\zeta))\cdot [\dot{f}_i (\zeta; a)
\bar{\dot{f}_j(\zeta; a)}] 
\end{equation}
$$+  r_{w_i,
w_j} (f_{0}(\zeta))\cdot [\dot{f}_i(\zeta; a) \dot{f}_j(\zeta; a)]\} d\theta-|\dot{f}_1(\zeta; a)|^2 < 0.$$
Note that the strict inequality follows because we cannot have
$\dot{f}(\zeta; a) \equiv 0, \zeta \in \partial \triangle$ by the maximum modulus principle, since $\hat{f}_{ta}(0)=(0,ta)$ so that $(\hat{f}_{ta})_j(0) = ta_j, \ j = 2,\ldots, n,$ and some $a_j  \neq 0$.

\end{proof}

\begin{corollary}
	\label{cor:3.3}
In the coordinates $w = (w_1,\ldots, w_n)$ near $\hat{{\mathcal K}}(0, 0) \in \CC\PP^n \setminus \bar{D}$ we have $$V_{\bar {D}} = \log \frac{1}{|w_1|} + R_D(w),$$ where $R_D(w)$ is smooth across $H_{\infty}$ which is locally given by $\{w_1 =
0\}$.
\end{corollary}

\begin{proof} 
This can be read off directly from the formula for $\hat{{\mathcal K}}$ near $(0, 0)$ in section 2, using the fact that $\hat{{\mathcal K}}$ is invertible across ${\mathcal E}$.
\end{proof}

We can interpret $H(0,\eta)$ using a ``Robin function'' associated to $G=G_{D'}$. Let $$r_G(z_1,...,z_n):=\limsup_{t\to 0} \bigl[G(tz_1,...,tz_n)-\log |t|\bigr].$$
Then $H(0,\eta)=r_G(1,\eta)$. Indeed, since $\log |\zeta|=G( f(\zeta; \eta))$ from (\ref{eqn:grkob}) we have
$$0=G( f(\zeta; \eta))-\log |\zeta|=G(\lambda_{(1,\eta)}(1, \eta)\zeta +0(|\zeta|^2))-\log |\zeta|$$
which shows, since $f(\zeta; \eta)$ is smooth in $\zeta$ and $\eta$,  that $r_G(\lambda_{(1,\eta)}(1, \eta))=0$ and the limit in the definition of $r_G$ exists. From the definition of $r_G$, it is a logarithmically homogeneous function, so that 
$$r_G(\lambda_{(1,\eta)}(1, \eta))=r_G((1, \eta))+\log \lambda_{(1,\eta)}=0$$
and $r_G((1, \eta))=-\log \lambda_{(1,\eta)}=H(0,\eta)$ from (\ref{ho}). Thus the Robin function $r_G$ coincides with $-\log \lambda$ at all points; hence we write $\lambda(v)$ where $v\in \CC^n$. Here we are identifying $T_0'(D')$ with $\CC^n$. We can interpret $r_G$ as the Green's function $G_I:= G_{I_0(D')}$ for $I_0(D')$. From Corollary 3.2 $I_0(D')$ is strictly convex; it is also balanced. If $\Omega$ is balanced and pseudoconvex; i.e., $\Omega=\{z\in \CC^n: u(z)<1\}$ where $u$ is psh and homogeneous: $u(tz)=|t|u(z)$, then $G_{\Omega}(z)=\log u(z)$. Since $I_0(D')=\{v:K_0(v)<1\}$, we have $G_I(v)=\log K_0(v) =-\log \lambda(v)=r_G(v)$ for $v\in I_0(D')$ as claimed. Moreover, the Kobayashi geodesics through the origin in $I_0(D')$ are flat disks. We remark in passing that (cf., \cite{lem3})
\begin{equation}
\label{kobe}
I_0(D')=\{v=f'(0):f \ \hbox{holomorphic on} \ \triangle, \ f(\triangle)\subset D', \ f(0)=0\}.
\end{equation}

We can use $G_I$ to define a Kelvin transformation $K_I$ from $I_0(D')$ to the complement of the closure of a balanced, strictly lineally convex domain $R(D)$. From the preceeding results, this map extends to $\hat K_I:\hat {I_0(D')} \to \CC\PP^n \setminus R(D)$. We will discuss the set $R(D)$ in the next section; to motivate this we recall the circular representation of the domain $D'$ (see \cite{lem3}). This is simply the linearization $\Phi=\Phi_{D'}$ from $\bar D'$ to $\bar {I_0(D')}$ described as follows: set $\Phi(0)=0$ and for $z\in \bar D'\setminus \{0\}$, let $f:\bar{\triangle}\to \bar D'$ be a Kobayashi geodesic in $D'$ with $f(0)=0$ and $f(t)=z$ for some $0<t\leq 1$. Let $g:\bar{\triangle}\to \bar {I_0(D')}$ be a Kobayashi geodesic in $I_0(D')$ with $g(0)=0$ and $g'(0)=\lambda f'(0)$ for some $\lambda >0$; define $\Phi(z)=g(t)$. Then $\Phi$ is a homeomorphism of $\bar D'$ onto $\bar {I_0(D')}$ which is smooth away from the origin. Extending $\Phi$ as a map $\hat \Phi$ on  the blow-up $\hat {\bar D'}$ of $\bar D'$ at the origin to the blow-up $\hat {\bar {I_0(D')}}$ of $\bar {I_0(D')}$ at the origin by requiring that $\hat \Phi$ fix the exceptional divisor gives an everywhere smooth map. 

For $v\in I_0(D')$, the Kobayashi geodesic through the origin in $D'$ in the direction of $v$ can be written as 
$$g_v(\zeta)=\lambda(v)v\zeta +0(|\zeta|^2)= e^{-r_G(v)}v\zeta +0(|\zeta|^2), \ \zeta \in \triangle;$$
the Kobayashi geodesic through the origin in $I_0(D')$ in the direction of $v$ can be written as 
$$\tilde g_v(\zeta)=\lambda(v)v\zeta= e^{-r_G(v)}v\zeta, \ \zeta \in \triangle.  $$
We can identify the tautological line bundle over $\PP(T_0(D'))$ (i.e., over $\CC\PP^{n-1}$) with the blow-up of $\CC^n$ at $0$; with this identification, $\zeta v\in \bar {I_0(D')}$ with $\zeta \in \triangle$ lifts to
$$\bigl(e^{-r_G(v)}\zeta(v_1,v_2,\ldots,v_n),[v_1:v_2:\ldots:v_n]\bigr)=(e^{-r_G(v)}\zeta v,[v])\in \hat {\bar {I_0(D')}}$$
and 
$$\hat \Phi (g_v(\zeta),[v])=(e^{-r_G(v)}\zeta v,[v])= (\tilde g_v(\zeta),[\frac{\tilde g_v(\zeta)}{\zeta}]).$$
In affine coordinates, we can consider $\partial I_0(D')/\sim$ as a parameter space for the Kobayashi geodesics, where $\sim$ denotes the equivalence relation from the circle action: $v\sim w$ if and only if $v=e^{i\theta}w$ for some $\theta$. For $v\in \partial I_0(D')$, we have $\lambda(v)=1$ so we can write the Kobayashi geodesic through the origin in the direction of $v$ as $g_v(\zeta)=v\zeta +0(|\zeta|^2)$. Then $\Phi(g_v(\zeta))=\zeta v$. The inverse map, $\Phi^{-1}(\zeta v)=g_v(\zeta)$, can be thought of as an exponential map from $\bar {I_0(D')}$ onto $\bar D'$.  

\vskip 5mm

\section{\bf Variational problem at $H_{\infty}$, and the Robin indicatrix.}
	\label{sec:varprob}
\vskip 3mm

Let $D, D'$ be as above. In \cite{blm} it was pointed out that the Monge-Amp\`ere solution $V_{\bar {D}}$ gave a  foliation of $\CC^n \setminus \bar{D}$ by holomorphic curves which could be represented as 
$h: \CC \setminus \bar{\triangle} \rightarrow \CC^n \setminus \bar{D}$, with Laurent expansion $$h(\zeta) =(h_1(\zeta),...,h_n(\zeta))= a_1 \zeta  + \sum_{j \leq 0} a_j \zeta^j, \ a_j\in \CC^n,$$ with
$a_1 \neq 0.$ It follows that the map extends holomorphically to a map, still denoted $h$, from $\CC\PP^1 \setminus \bar{\triangle}$ to $\CC\PP^n \setminus
\bar{D}$ with $h(\infty) = [0:a_1] \in H_{\infty} \simeq (\CC^n \setminus \{0\})/\CC^{*}.$ Replacing $\zeta$ by $\frac{1}{\zeta}$, we can consider $h$ as a holomorphic map $f$ 
from $\triangle \rightarrow \CC\PP^n \setminus \bar{D}$, with $f(0) = [0:a_1] \in H_{\infty}.$ Without loss of generality, we will assume $a_1 = (a_{1,1},\ldots,
a_{1,n})$, with $a_{1,1} \neq 0$, so that we can write this out in terms of the local coordinates $w = (w_1,\ldots, w_n)$ (see section 2) near $[0:a_1], \ f(\zeta) =
(f_1(\zeta),\ldots, f_n(\zeta))$ with $f(0) = (1, \frac{a_{1,2}}{a_{1,1}},\ldots, \frac{a_{1,n}}{a_{1,1}})$ and 
$$f_1(\zeta)=\frac{1}{h_1(1/\zeta)} =\frac{\zeta}{a_{1,1}+\sum_ {j\leq 0}a_{j,1}\zeta^{1-j}}$$
so that $f_1'(0)=\frac{1}{a_{1,1}}$.

Now, given a point $c = [(1, c_2,\ldots, c_n)] \in H_{\infty}$ consider the class $\mathcal{H}_{(c, D)}(\triangle)$ of all  holomorphic disks 
$f: \triangle \rightarrow \CC\PP^n$ which are in $\mathcal{C}^{k, \alpha}(\bar{\triangle})$, $k, \alpha$ as before, and are such that
\begin{equation}
\label{eqn:defncompetitors}
\begin{array}{rl}
i) & f(0) = c \\
	&	\\
	ii) & f(\partial \triangle) \subset \bar {D} \\
	&	\\
	iii) & f_1'(0) = 1/\rho > 0 \ \hbox{in} \ w \ \hbox{coordinates}.
\end{array}
\end{equation}

\vskip 2mm

\noindent Note that $\mathcal{H}_{(c, D)}(\triangle)\not =\emptyset$ since there clearly exist maps $f$ satisfying i)-iii) for $\rho$ sufficiently small.

By Theorem \ref{th:kelinf}, we know that there is a unique foliation disk $f = f_c$ passing through $c$ which is in the class 
$\mathcal{H}_{(c, D)}(\triangle)$. For any $f \in \mathcal{H}_{(c, D)}(\triangle)$ consider the function $f^{*}V_{\bar {D}} + \log |\zeta|$. Note that it is subharmonic
and continuous on all of $\triangle$, $\mathcal{C}^{k, \alpha}$ at $\partial \triangle$  and is $\leq 0$ on $\partial \triangle$. Set $$\mu(f) = \,\mbox{value of}
\, f^{*}V_{\bar {D}} + \log |\zeta| \, \mbox{at} \, \zeta = 0.$$ Notice that we have $\mu(f) \leq 0$, with equality holding if and only if $f = f_c$. Explicitly, from Corollary 3.3, 
\begin{equation}
\label{eqn:muf}
\mu(f)
= \log \rho(f) + R_D(0, c_2,\ldots, c_n),
\end{equation} 
so that $f$ maximizes $\mu(f)$ if and only if it maximizes $\rho(f)$ among all competitors. We set 
$$\mu(c) = \mu(c, D) = \mu(f_c) = \max \{\mu(f) | f \in \mathcal{H}_{(c, D)}(\triangle)\} <
0,$$ 
and equivalently 
$$\rho(c) = \rho(c, D) = \rho(f_c) = \max \{\rho(f) | f \in \mathcal{H}_{(c, D)}(\triangle)\} > 0.$$

Finally, one does not need to know $V_{\bar D}$ to express the variational problem. If we take $w$ to be any local holomorphic defining function for $H_{\infty}$, then maximizing $\mu(f)$ is equivalent to maximizing 
\begin{equation}
\label{eqn:mutilde}
\tilde{\mu}(f) = 
\,\mbox{value of}
\, \bigl[\log(\frac{1}{|f^*w|}) - \log(\frac{1}{|\zeta|})\bigr] \, \mbox{at} \, \zeta = 0.
\end{equation}
over $f \in \mathcal{H}_{(c,D)}$.

We would like to package this information into an exponential map for our functional, analogous to the inverse of the map $\Phi_{D'}$ in the previous section. First, recall the definition of the {\it Robin function} $\rho_{\bar D}$ (see \cite {sic}): for $z=(z_1,...,z_n)\in \CC^n$, 
$$\rho_{\bar D}(z_1,...,z_n)=\limsup_{t\to \infty} \bigl[ V_{\bar D}(tz_1,...,tz_n)-\log |t|\bigr].$$
For $D$ as above, the limit exists and defines a logarithmically homogeneous psh function: $\rho_{\bar D}(tz)=\rho_{\bar D}(z) +\log |t|$. 

Let $c=[0:1:c_2:\ldots:c_n]\in H_{\infty}$, and let $f_c$ be the unique foliation disk through $c$ so that the conditions of (\ref{eqn:defncompetitors}) are satisfied for maximal $\rho=\rho(f_c)$.  Then from (\ref{eqn:muf})
\begin{equation}\label{rhorobin}\rho(f_c)=\rho(c,D)=e^{-R_D(0,c_1,c_2\ldots,c_n)}=e^{-\rho_{\bar D}(1,c_1,\ldots,c_n)};\end{equation}
here the domain of $R_D$ is given in local $w$-coordinates and the domain of $\rho_{\bar D}$ is given in affine $z$-coordinates.
\par
In affine coordinates in $\CC^n=\CC\PP^n\setminus H_{\infty}$, a foliation disk $f$ satisfies the growth condition $\zeta f(\zeta)=O(1)$, hence is of the form
\begin{equation}\label{extrdisk} f(\zeta)=\frac{a}{\zeta} + \sum_{j=0}^{\infty}b_j\zeta^j, \ \ \zeta\in\Delta\setminus\{0\},\end{equation}
where $a=(a_1,\ldots,a_n)$ and $b_j=(b_{j,1},\ldots,b_{j,n})$ for all $j$; for such an extremal we have $f(0)=[0:a_1:\ldots:a_n]=[0:a]\in H_{\infty}$.  Uniqueness of the foliation disk implies that modulo the circle action, $a$ is \emph{unique}.  If in fact $[0:a]=c$ with $c$ as above, then $(\ref{eqn:defncompetitors}({\it iii}))$ implies that in affine coordinates on $\CC^n=\CC\PP^n\setminus H_{\infty}$, the  first component of $f_c$ has the form $[f_c(\zeta)]_1=\frac{\rho(c,D)}{\zeta} +O(1)$.
If $f$ is parametrized so that $a_1\in\RR_+$, then by (\ref{extrdisk}) and the fact that $f=f_c$  we see that $\rho(c,D)=a_1$.  It follows by homogeneity that $\rho(c,D)c_j=a_j$ for all $j$, so that (\ref{extrdisk}) can be rewritten as
\begin{equation}\label{extrdisk1}
f_c(\zeta) = f_{(c,D)}(\zeta)= \frac{\rho(c,D)(1,c_2,\ldots,c_n)}{\zeta} + \sum_{j=0}^{\infty}b_j\zeta^j.
\end{equation}

More generally, for any $v\in\CC^n\setminus \{0\}$, consider $c_v:=[0:v]\in H_{\infty}$, and let $f_{c_v}$ be the  extremal disk through $c_v$ given in the form
\begin{equation}\label{eqn:rhov}
f_{c_v}(\zeta)=\frac{\tilde \rho(v)v}{\zeta} + \sum_{j=0}^{\infty}b_j\zeta^j
\end{equation}
for some appropriate $\tilde \rho(v)>0$. From (\ref{eqn:grkob}) and (\ref{eqn:greensiciak}) (see \cite{lem2}), $V_{\bar D}(f_{c_v}(\zeta))=\log\frac{1}{|\zeta|}$ and we have
$$0=V_{\bar D}(f_{c_v}(\zeta))-\log\frac{1}{|\zeta|}=V_{\bar D}\bigl(\frac{\tilde \rho(v)v}{\zeta} + \sum_{j=0}^{\infty}b_j\zeta^{j}\bigr)-\log\frac{1}{|\zeta|}$$
so that $\rho_{\bar D}(\tilde \rho(v)v)=0$. By logarithmic homogeneity of $\rho_{\bar D}$, we have 
$$\rho_{\bar D}(\tilde \rho(v)v)=\rho_{\bar D}(v)+\log \tilde \rho(v)=0$$
so that 
\begin{equation}\label{rhov}
\tilde \rho(v)=\exp{[- \rho_{\bar D}(v)]}.
\end{equation}
Observe that equation (\ref{rhorobin}) is just a special case of this equality.

We will soon see that the set $R(D)$ from the previous section can be identified as
\begin{equation}\label{RDrho}R(D)=\{z\in \CC^n: \rho_{\bar D}(z)\leq 0\}.\end{equation}

Define a map $F_D$ from a neighborhood ${\mathcal O}_D$ of $H_{\infty}\subset\CC\PP^n$ to $\CC\PP^n\setminus\bar D$ by the equation
\begin{equation}
\label{eqn:robexp}
F_D\bigl([\frac{\zeta}{e^{-\rho_{\bar D}(v)}}:v]\bigr)=f_{c_v}(\zeta),\ \ |\zeta|<1\, .
\end{equation}
Note that $F_D([0:v])= [0:v]$; i.e., $F_D$ is well-defined as the identity map on $H_{\infty}$.  We give the verification that (\ref{eqn:robexp}) indeed gives a well-defined mapping 
$F_{D}: \mathcal{O}_D \rightarrow \CC\PP^n \setminus \bar{D}$ 
in Remark 2 at the end of this section.

When $\zeta\neq 0$, $[\frac{\zeta}{e^{-\rho_{\bar D}(v)}}:v]$ is given in affine coordinates by $z=\frac{e^{-\rho_{\bar D}(v)}v}{\zeta}$.  Now
$$
\rho_{\bar D}(z)=\rho_{\bar D}(v)+\log\frac{e^{-\rho_{\bar D}(v)}}{|\zeta|} = \log\frac{1}{|\zeta|},
$$
which, since $|\zeta|<1$, shows that ${\mathcal O}_D\cap \CC^n=\{\rho_{\bar D}>0\}=\CC^n\setminus R(D)$, provided (\ref{RDrho}) holds.

Note that for $c=[0:1:c_2:...:c_n]\in H_{\infty}$, 
\begin{equation}
\label{eqn:cs} 
[\zeta/\rho(c,D):1:c_2:...:c_n]=[1/\rho(c,D):1/\zeta:c_2/\zeta:...:c_n/\zeta].
\end{equation}
We will abuse notation and write (\ref{eqn:cs}) as $c/\zeta$. Thus we write, using (\ref{extrdisk1}), 
\begin{equation}
\label{eqn:RhoC} 
F_D(c/\zeta):=F_D([\zeta/\rho(c,D):c])=f_c(\zeta) = \frac{\rho(c,D)(1,c_2,\ldots,c_n)}{\zeta} + \sum_{j=0}^{\infty}b_j\zeta^j.
\end{equation}

We provide justification for this notational abuse as follows. As with the circular representation $\Phi$ mapping $D'$ onto its linearization $I_0(D')$ in the previous section, we can use affine coordinates to interpret $F_D$ as an exponential map.  We call ${\mathcal O}_D=\CC\PP^n\setminus R(D)$ the \emph{Robin indicatrix}, and $F_D$ the \emph{Robin exponential map}. The set $\partial R(D)/\sim$, where again $\sim$ denotes the equivalence relation from the circle action, can be used as a parameter space for the leaves of the foliation of $\CC^n \setminus \bar D$, and from (\ref{rhov}) and (\ref{RDrho}) $\tilde \rho(v)=1$ for $v\in \partial R(D)$. For such $v$, and $0<|\zeta|<1$, the point $v/\zeta\in \CC^n \setminus R(D)$. Thus from (\ref{eqn:rhov}) we can consider $F_D:\CC^n \setminus R(D)\to \CC^n \setminus \bar{D}$ via 
\begin{equation}
\label{eqn:Rhov} 
F_D(v/\zeta) = v/\zeta +  \sum_{j=0}^{\infty}b_j\zeta^j;
\end{equation} 
i.e., $F_D^{-1}(v/\zeta +  \sum_{j=0}^{\infty}b_j\zeta^j)=v/\zeta$ is the affine linearization from $\CC^n \setminus \bar{D}$ onto $\CC^n \setminus R(D)$. 

We now verify (\ref{RDrho}).  Recall from the previous section that one can also construct a Kelvin transform $K_I$ defined on the Kobayashi indicatrix $D_0':=I_0(D')$, considered as a subset of $\CC^n$, the dual of $D_0=R(D)^o$.  We then get an extended Kelvin transform $\hat K_I:\hat D_0'\longrightarrow \CC\PP^n\setminus\bar D_0$. Since  $D_0'$ is circled, it follows easily that $D_0$ is circled.

Let $G_I$ be the Green function of $D_0'$ with logarithmic pole at the origin, so that, from (\ref{eqn:greensiciak}), $V_{\bar D_0}(\hat K_I(z))=-G_I(z)$. We also need the following facts:
\begin{itemize}
\item[(i)] $G_I(z)=r_G(z)=\lim_{\zeta\to 0} [G(\zeta z) -\log|\zeta|]$.
\item[(ii)] $\hat K_I(z)=\lim_{\zeta\to 0} \zeta {\mathcal K}(\zeta z)$.
\end{itemize}

Formula (i) was proved in the previous section. For (ii), define the function $h(\zeta,z):=G_{D'}(\zeta z)-\log|\zeta|=G(\zeta z)-\log|\zeta|$.  Then $h$ is smooth whenever $\zeta$ and $z=(z_1,...,z_n)$ are nonzero.  When $z_1\neq 0$, we have, from (\ref{eqn:pistar}), 
$$
h(\zeta,z)=H(\zeta z_1,\eta)+\log|z_1|, \ (\eta_j=z_j/z_1, \ j=2,...,n). 
$$
Moreover, $H$ extends smoothly across the exceptional divisor; i.e., $h$ extends smoothly across $\zeta=0$. Then $G_I(z)=h(0,z)$, and we can differentiate inside the limit to obtain
\begin{eqnarray*}
\frac{\partial}{\partial z_j}G_I(z) = \frac{\partial h}{\partial z_j}(0,z)
&=&\lim_{\zeta\to 0} \frac{\partial h}{\partial z_j}(\zeta,z) \\
&=&\lim_{\zeta\to 0} \frac{\partial}{\partial z_j}\bigl(G(\zeta z)-\log|\zeta|\bigr)
=\lim_{\zeta\to 0} \frac{\partial G}{\partial z_j}(\zeta z)\cdot\zeta.
\end{eqnarray*}
Plugging this last expression into the formula for $\hat K_I$ (cf., equation (\ref{Kelvin})) yields (ii).

Using (i), (ii), (2.2) and the continuity of $\rho_{\bar D}$ we have
\begin{eqnarray*}
\rho_{\bar D}(\hat K_I(z))=\lim_{\zeta\to 0} \rho_{\bar D}(\zeta {\mathcal K}(\zeta z))
&=& \lim_{\zeta\to 0} [\rho_{\bar D}({\mathcal K}(\zeta z))+\log|\zeta|]   \\
&=& \lim_{\zeta\to 0} [V_{\bar D}({\mathcal K}(\zeta z))+\log|\zeta|] \\
&=& \lim_{\zeta\to 0} [-G(\zeta z)+\log|\zeta|] = -G_I(z).
\end{eqnarray*}
But $-G_I(z)=V_{\bar D_0}(\hat K_I(z))=\rho_{\bar D_0}(\hat K_I(z))$ on $\PP^n\setminus D_0$. Hence  $\rho_{\bar D}=\rho_{\bar D_0}$, which proves (\ref{RDrho}).

Analogous to (\ref{kobe}), we have
$$R(D)=\{v\in \CC^n: v=\lim_{\zeta \to 0}\zeta f(\zeta), \ f(\partial \triangle) \subset \bar D\}.$$

\vskip 3mm

\noindent {\bf Remark 1.} This functional (\ref{eqn:mutilde}) and the stability properties of its second variation, analogous to (\ref{eqn:2varsign}) above, appear already in \cite{dbx}.

\noindent {\bf Remark 2.} We verify that (\ref{eqn:robexp}) is a well-defined mapping $F_{D}: \mathcal{O}_D \rightarrow \CC\PP^n \setminus \bar{D}$.

Let $v,\tilde v\in\CC^n$, $v=\lambda\tilde v$.  Then $c_v=c_{\tilde v}$ so that $f_{c_v}$ and $f_{c_{\tilde v}}$ defined by (\ref{eqn:rhov}) parametrize the same extremal disk.  Hence $f_{c_{\tilde v}}(\zeta)=f_{c_v}(e^{i\theta}\zeta)$ for some $\theta=\theta(\lambda)$.  From the first term in (\ref{eqn:rhov}), we have
$\frac{\rho(v)v}{e^{i\theta}\zeta}=\frac{\rho(\tilde v)\tilde v}{\zeta}$, which implies that $\rho(v)=\frac{\rho(\tilde v)}{|\lambda|}$ and $\theta(\lambda)=\arg(\lambda)$.  Now
$$
[\frac{\zeta}{e^{-\rho_{\bar D}(\tilde v)}}:\tilde v] = [\frac{\zeta}{\rho(\tilde v)}:\tilde v]
= [\frac{\zeta}{|\lambda|\rho(v)}:\tilde v] = [\frac{e^{i\theta(\lambda)}\zeta}{\rho(v)}:v] 
=[\frac{e^{i\theta}\zeta}{e^{-\rho_{\bar D}(v)}}:v]\, ;
$$
thus (\ref{eqn:robexp}) is well-defined.

We also have a bundle interpretation of the map $F_D$.  Consider the normal bundle of $H_{\infty}$:

$$N_{H_{\infty}}:=\{(c,\nu): c\in H_{\infty}, \ \nu\in (T_cH_{\infty})^{\perp}\}.$$
Suppose we are in a region of $H_{\infty}$ where $w$-coordinates are valid; here $(0,c_2,\ldots,c_n)$ corresponds to $[0:1:c_2:\ldots:c_n]$.  Local coordinates for $N_{H_{\infty}}$ may be given by $(\zeta,c_2,\ldots,c_n)$ where $\zeta$ is the fiber variable.  We can identify ${\mathcal O}_D$ from (\ref{eqn:robexp}) with a subset of $N_{H_{\infty}}$ via
$$
N_{H_{\infty}}\ni (\zeta,c_2,\ldots,c_n)\sim [\frac{\zeta}{e^{-\rho_{\bar D}(1,c_2,\ldots,c_n)}}:1:c_2:\ldots:c_n]\in {\mathcal O}_D.
$$
The same calculation used in showing (\ref{eqn:robexp}) is well-defined shows that different local trivializations in a neighborhood of a base point $[0:v]\in H_{\infty}$ give the same point in ${\mathcal O}_D$ under this identification. We can now reinterpret $F_D$ as a map from a subset of $N_{H_{\infty}}$ to $\CC\PP^n\setminus D$.

 \vskip 5mm

\section{\bf Passing to a real convex body as limit.}
	\label{sec:convexlimit}

\vskip 3mm

Let $K \subset \subset \RR^n$ be a compact convex body, and let $V_K$ be the Siciak-Zaharjuta extremal function on $\CC^n$. In
\cite{blm} it  was shown that through every point $z \in \CC^n \setminus K$ there passes a holomorphic curve $f(\triangle)$ as in the previous section such that $f^{*}V_K = \log
\frac{1}{|\zeta|}$. This information was derived using a decreasing sequence of strictly convex $\sigma$-invariant open sets $D_j$ such that $\cap D_j = K$, and a
normal families argument on a sequence of foliation curves $f_j(\triangle)$ for $\CC^n \setminus D_j$ (recall $\sigma$ is the usual complex conjugation of $\CC^n$.) It
was also shown that all such curves on $\CC \setminus \triangle$ were of the form
$$h: \CC \setminus \triangle \rightarrow \CC^n \setminus K,$$
$$h(\zeta) = a_1 \zeta + a_0 + \bar{a}_1 \frac{1}{\zeta}, \ a_0\in \RR^n, \ a_1\in \CC^n\setminus \{0\}.$$ As noted in section 2 already, this curve extends through the point $[0:a_1] \in
H_{\infty}$,  for $\zeta = \infty$. We would like to extend the variational properties of section 4 to these curves.

We saw in Corollary \ref{cor:3.3} that for $D$ as in the previous two sections, $V_{\bar {D}} = \log \frac{1}{|w_1|} + R_D(w)$, where $R_D(w)$ is smooth across $H_{\infty}$. This local coordinate expression
is equivalent to saying  that the function $V_{\bar {D}}(z) - \log |z|$ on $\CC^n \setminus \{0\}$ extends smoothly across $H_{\infty}$. The smoothness strengthens a result of Siciak \cite {sic} in this  special case: Siciak's result is that the Robin function $\rho_K$ associated to a compact set $K\subset \CC^n$ with $V_K$ continuous is itself continuous. We recover this fact for a convex body $K\subset \RR^n$.

\begin{theorem}
	\label{th:robincont}
	For $K$ a convex body in $\RR^n$, $V_K(z) - \log |z|$ extends continuously across $H_{\infty} \subset \CC\PP^n$.
\end{theorem}

\begin{proof} Let $D_j$ be a decreasing sequence of bounded, smoothly bounded strictly convex open sets in $\CC^n$ such that $\cap D_j = K$. Then it is well 
known that the extremal functions $V_{\bar {D}_j}$ are continuous, monotonically increasing and converge uniformly to $V_K$. Then $V_{\bar {D}_j}(z) - \log |z|$
converge uniformly to $V_K(z) - \log |z|$ on $\CC^n \setminus \{0\},$ and by density of $\CC^n$ in $\CC\PP^n$, the extensions of the functions $V_{\bar {D}_j}(z) - \log |z|$ converge 
monotonically and uniformly on $\CC\PP^n \setminus \{0\}$. Thus, the function $V_K(z) - \log |z|$ has a continuous extension across $H_{\infty}$. \end{proof}

Fixing a point $c \in H_{\infty}$, which for convenience we assume is in the domain of the coordinates $w = (w_1,\ldots, w_n)$, we can consider the space 
$\mathcal{H}_{(c, K)}(\triangle)$ of all holomorphic maps $f: \triangle \rightarrow \CC\PP^n$ with $f(\partial
\triangle)
\subset K$ which are of the form (in affine $z-$coordinates on $\CC^n$) $$f(\zeta) = h(1/\zeta) = \frac{\rho}{\zeta} (1,
c_2,\ldots, c_n) + \rho \zeta (1, \bar{c}_2,\ldots, \bar{c}_n) + a_0,$$ where $\rho > 0$ and $a_0 \in \RR^n$. Note in the coordinates $w$ around $c \in
H_{\infty}$,

\vskip 2mm

\centerline{$\begin{array}{rclc}
w_1 & = f_1(\zeta)=\frac {1} {h_1(1/\zeta) }& =   \zeta \cdot \frac{1}{\rho+ a_{0,1}  \zeta + \rho \zeta^2}, &\\
	&	&	\\
w_j 	& = f_j(\zeta)=\frac {h_j(1/\zeta)}{h_1(1/\zeta)}& =  \frac{\rho \bar c_j + a_{0, j} \zeta + \rho c_j \zeta^2}{\rho + a_{0,1}  \zeta + \rho \zeta^2}, \; j = 2,\ldots, n.& 
\end{array}$}

\vskip 2mm

This is a finite dimensional set of mappings, and for such maps, $f(0) = c$ and $ f'_1(0) = 1/\rho > 0$ in the $w$ coordinates. 
Theorem \ref{th:robincont} states that $V_K=\log \frac{1}{|w_1|} + R_K(w)$ where $R_K(w)$ is continuous across $H_{\infty}$. The variational interpretation carries over
to this limit case, that is:

\begin{theorem}
	\label{th:5.2}
	The mapping $f \in \mathcal{H}_{(c, K)}(\triangle)$ satisfies $$f^{*}V_K = \log \frac{1}{|\zeta|}$$
	if and only if $\rho = \rho(f) > 0$ is maximal for maps in $ \mathcal{H}_{(c, K)}(\triangle)$.
\end{theorem}

\begin{proof} For any $f \in  \mathcal{H}_{(c, K)}(\triangle)$, we can write $$f^{*}V_K = \log \frac{1}{|\zeta|} + \log \rho(f) + R_K(f(\zeta)) + O(|\zeta|).$$ 
The function $f^{*}V_K - \log \frac{1}{|\zeta|}$ is continuous and subharmonic on $\triangle$ and $\leq 0$, so we conclude
$$\log \rho(f) \leq -R_K(f(0)),$$ with equality if and only if $f^{*}V_K = \log \frac{1}{|\zeta|}.$
\end{proof}

	Note that the extremal value $\rho(f) = \rho(c, K) = e^{- R_K(0,c)}=e^{-\rho_K(1,c_2,...,c_n)}$ is the limit of the extreme values $\rho(c, D_j) = e^{- R_{D_j}(0,c)}=e^{-\rho_{\bar D_j}(1,c_2,...,c_n)}$ 
for any sequence $D_j \rightarrow K$, and that this limit is uniform in $c \in H_{\infty}$. In affine coordinates, the Robin functions $\rho_{\bar D_j}$ converge uniformly to the Robin function $\rho_K$. We have the diffeomorphisms $F_{D_j}: \mathcal{O}_{D_j} \to \CC\PP^n \setminus \bar D_j$, and similarly we can define $F_K: \mathcal{O}_K = \cup \mathcal{O}_{D_j} \to \CC\PP^n \setminus K$. Define 
\begin{equation}
\label{krho}
K_{\rho}:=\{z\in \CC^n: \rho_K(z)\leq 0\}. 
\end{equation}
Then the sets $R(D_j)$ decrease to $K_{\rho}$ and we consider $F_K:\CC\PP^n \setminus K_{\rho} \to \CC\PP^n \setminus K$. As in 
(\ref{extrdisk1}), for $c=[0:1:c_2:...:c_n]\in H_{\infty}$ we write $f_c=f_{(c,K)}$ where
\begin{equation}
\label{eqn:fc}
\begin{array}{rcl}
f_c(\zeta) & = & f_{(c,K)}(\zeta) \\
	& = & \rho(c, K) \{(1, c_2,\ldots, c_n) \zeta^{-1} + (1, \bar{c}_2,\ldots, \bar{c}_n) \zeta\} + a_0(c, K).
\end{array}
\end{equation}
Then as in (\ref{eqn:RhoC}) we use the notation
\begin{equation}
\label{eqn:fkc}
F_K(c/\zeta )=\rho(c, K) \{(1, c_2,\ldots, c_n) \zeta^{-1} + (1, \bar{c}_2,\ldots, \bar{c}_n) \zeta\} + a_0(c, K).
\end{equation}
Analogous to (\ref{eqn:Rhov}), if we consider $\partial K_{\rho}$ modulo the circle action as a parameter space for our extremal curves, since $\rho_K(v)=0$ for $v\in \partial K_{\rho}$ from (\ref{krho}), we can consider $F_K: \CC^n \setminus K_{\rho} \to \CC^n \setminus K$ via
\begin{equation}
\label{eqn:fk}
F_K(v/\zeta)=a_0(v,K)+ \bigl(v/ \zeta + \bar v\zeta\bigr).
\end{equation}
The next remark shows that convex bodies $K\subset \RR^n$ are natural sets to consider.

\vskip 2mm

\noindent {\bf Remark}. Let $\{D_j\}$ be a decreasing sequence of relatively compact, strictly lineally convex domains in $\CC^n$ that are invariant under conjugation (i.e., $\sigma(D_j)=D_j$ for all $j$), and suppose
$K:=\bigcap_j D_j$ is a compact set contained in $\RR^n$.  Then $K$ is convex.

To see this, for each $j$, we let $K_j:=\bar D_j\cap\RR^n$.  Then $K_j\downarrow K$, so to show convexity of $K$ it suffices to show that $K_j$ is convex for each $j$.  To this end, fix $j$ and let $a\in\partial K_j$.  
We may assume $a=0$.  Then
 $$T_0^{\CC}(\partial D_j)\cap \bar D_j=\{0\},$$
and we can write 
$$T_0^{\CC}(\partial D_j)=\{z=(z_1,\ldots,z_N)\in\CC^n: \sum_{k=1}^N b_kz_k=0\}$$
for some $b_k\in {\CC}$. On the other hand, by symmetry, we have
$$T_0^{\CC}(\partial \sigma(D_j))=\{z=(z_1,\ldots,z_N)\in\CC^n: \sum_{k=1}^N \bar b_kz_k=0\}.$$
Thus if $D_j=\sigma(D_j)$, then we can take $b_k=\bar b_k$ for all $k$.  It follows that $T_0^{\CC}(\partial D_j)\cap\RR^n$ is a real hyperplane whose intersection with $K_j$ is $\{0\}$.  Since $0$ was an arbitrary point of $\partial K_j$, it follows that $K_j$ is convex.  Since $j$ was also arbitrary, we conclude that $K$ is convex. This argument only uses the existence, through each boundary point of $D_j$, of a complex hyperplane that does not intersect $D_j$; i.e., it is valid if we assume each $D_j$ is conjugation invariant and {\it weakly lineally convex}.

\vskip 5mm

\section{\bf Geometric interpretation of extrema.}
\label{sec:geomint}

\vskip 3mm

As already noted in \cite{blm}, the $f \in  \mathcal{H}_{(c, K)}(\triangle)$ extend by reflection to mappings of the Riemann sphere $\CC\PP^1 \rightarrow \CC\PP^n$. 
The image curves are algebraic, of degree two in $\CC\PP^n$. Let $\sigma$ be the usual conjugation $(z_1,\ldots, z_n) \rightarrow (\bar{z}_1,\ldots, \bar{z}_n)$ of
$\CC^n$ which we will also consider on $\CC\PP^n$, by extension. 

If $c \in H_{\infty}$ and $c = \sigma(c)$, then the competitor mappings $f \in \mathcal{H}_{(c, K)}(\triangle)$ are of the form 
$$f(\zeta) = \rho(1, c_2,\ldots, c_n)(\zeta + \zeta^{-1}) + a_0,$$ in affine coordinates, where we again assume we have normalized the coordinate $c_1 =
1$. In this case, $f: \CC\PP^1 \rightarrow \CC\PP^n$ double covers the complex projective line in $\CC\PP^n$ in the direction $(1, c_2,\ldots, c_n) \in \RR^n$ and
through the point $a_0 \in \RR^n.$ In particular, $f(\partial \triangle)$ is the real line segment $a_0 + t \cdot (1, c_2,\ldots, c_n), \ t \in [-2\rho, 2\rho]$. In particular, we see that we maximize $\rho(f)$ among competitors if and only if we maximize the length of the segment.
Thus, given $c \in H_{\infty}$ such that $\sigma(c) = c$, the extremal $\rho$ comes from parametrizing the maximal line segment in the direction $(1, c_2,\ldots,
c_n) \in \RR^n$ contained in $K$; we are free to vary $a_0$ and $\rho$ to achieve this maximization. The point $a_0$ is the center of the maximal segment and the segment's length is
$4\rho |(1, c_2,\ldots, c_n)|.$

If $\sigma(c) \neq c$, then for $f \in  \mathcal{H}_{(c, K)}(\triangle)$, the extension $f: \CC\PP^1 \rightarrow \CC\PP^n$ is a non-singular quadric curve which 
intersects $H_{\infty}$ at $c$ and $\sigma(c)$. If we write $f$ in affine coordinates, $f(\zeta) = \frac{\rho}{ \zeta} (1, c_2,\ldots, c_n) +  \rho \zeta
(1,\bar{c}_2,\ldots, \bar{c}_n) + a_0$, then $f(\CC\PP^1)$ lies in the projective closure of the affine complex plane parametrized as $a_0 + \rho [z 
\,\Re(1, c_2,\ldots, c_n) + w \,\Im(1, c_2,\ldots, c_n)], \ z, w \in \CC,$ whose intersection with $\RR^n$ is the real affine plane parametrized as $a_0 +
\rho [s \,\Re(1, c_2,\ldots, c_n) + t
\,\Im(1, c_2,\ldots, c_n)], \ s, t \in \RR.$ The intersection of $f(\CC\PP^1)$ with this real affine plane is the ellipse given parametrically by
$$f(e^{i\theta}) = a_0 + 2\rho \Re(e^{-i\theta} (1, c_2,\ldots, c_n)).$$ By applying a real orthogonal transformation to all of $\RR^n
\subset \CC^n \subset \CC\PP^n$, we can assume $c_3 = \ldots = c_n = 0,$ and $\Re(1, c_2,\ldots, c_n) = (1, \alpha, 0,\ldots, 0)$ and $\Im(1, c_2,\ldots, c_n) = (0,
\beta, 0,\ldots, 0)$, where $\alpha, \beta \in \RR$. Our ellipse is now the translate by $a_0 \in \RR^n$ of the planar ellipse in $\RR^2 \equiv \RR^2 \times \{0\}
\subset \RR^n$ given by $(x, y) = 2\rho (\cos \theta, \alpha \cos \theta + \beta \sin \theta)$. Calculating, the area of the ellipse bounded by this curve
is given by $$|\int_0^{2\pi} x(\theta) \dot{y}(\theta) \, d\theta | = 4\pi |\beta| \rho^2 = 4\pi |\Im(c)|\rho^2.$$ Notice that this last expression
is invariant under the orthogonal transformation we used to simplify the coordinates. To summarize, our parameters $c = (1, c_2,\ldots, c_n)$ determine a real
2-plane (spanned by $\Re(1, c_2,\ldots, c_n)$ and $\Im(1, c_2,\ldots, c_n)$), and the family of all ellipses with given directions for the major and minor axes and
the eccentricity. Among these ellipses, we seek to adjust the center $a_0$ of the ellipse and the ``scale factor" $\rho$ to maximize the area of the
ellipse among all such ellipses which are also contained in $K$.

\vskip 5mm

\section{\bf Uniqueness of extremal curves.}
\label{sec:uniqueness}

\vskip 3mm

It is well known that for ``degenerate" convex bodies $K \subset \RR^n$, there may be many holomorphic disks $f(\triangle)$ through a given point 
$z \in \CC^n \setminus K$ such that $f^{*}V_K = \log \frac{1}{|\zeta|}.$ It is now easy to describe when there is more than one extremal curve passing through
the same point $c \in H_{\infty}$. Let $f^{(i)}=f_c^{(i)}, \  i = 0, 1,$ be two such disks. From (\ref{eqn:fc}) these maps are given in affine coordinates by 
$$f^{(i)}(\zeta) = \frac{\rho(c, K)} {\zeta}(1, c_2,\ldots, c_n) + \zeta \rho(c, K)  (1, \bar{c}_2,\ldots, \bar{c}_n) + a^{(i)}_0, \ i = 0, 1.$$ Setting 
$$f^{(t)}(\zeta) = (1-t) f^{(0)}(\zeta) + t f^{(1)}(\zeta), \ t \in [0, 1],$$ one obviously has $f^{(t)}(0) = c, \ \rho(f^{(t)}) = \rho(c, K)$, and since
$f^{(i)}(e^{i \theta}) \in K, \ \theta \in [0, 2\pi], \ i = 0, 1,$ we conclude, by the convexity of $K$, that $f^{(t)}(e^{i\theta})$ lies in $K$ for all
$\theta \in [0, 2\pi]$. Thus, $f^{(t)} \in  \mathcal{H}_{(c, K)}(\triangle)$ for all $t \in [0, 1]$, and we have shown that the set of all extremals in $
\mathcal{H}_{(c, K)}(\triangle)$ is a convex set. Note that $$f^{(t)}(\zeta) = f^{(0)}(\zeta) + t (a^{(1)}_0 - a^{(0)}_0).$$ 
In particular, the set of centers $\{a_0(c,K)\}$ of extremal curves associated to $c\in H_{\infty}$ is a closed convex set.

We call a set $X\subset\partial K$ a {\it face} (or, if $K\subset\RR^2$, an {\it edge}) of $K$ if $X$ is not a singleton and there is an affine hyperplane $H$ such that $K\cap H=\partial K\cap H = X$. It is easy to see that a face of a convex body must also be convex. Note also that in this definition we allow faces of codimension greater than one. We call two faces {\it parallel} if they lie on parallel hyperplanes. We conclude the following result:

\begin{theorem}
	\label{th:uniqueness}
	If $\partial K$ does not contain two parallel faces, then for every $c \in H_{\infty}$, there is a unique extremal curve 
$f_c \in  \mathcal{H}_{(c, K)}(\triangle)$. In particular, if $K\subset \subset \RR^n$ is the closure of a smooth, strictly convex domain, then there is a
unique $f_c$ through every $c \in H_{\infty}$.
	\end{theorem}

Without loss of generality, we may assume, via translation, that the origin is an interior point of $K$. We can write $K=\cap_j D_j $ where $D_j$ is a nested sequence of bounded, smoothly bounded strictly lineally convex domains containing the origin and we can define $F_K$ as in (\ref{eqn:fkc}). Thus $F_K: \CC\PP^n\setminus K_{\rho} \to \CC\PP^n \setminus K$ and we write
$$F_K(c/\zeta )=\rho(c, K) \{(1, c_2,\ldots, c_n) \zeta^{-1} + (1, \bar{c}_2,\ldots, \bar{c}_n) \zeta\} + a_0(c, K).$$

\begin{corollary}
	\label{cor:uniqueness}
	With $K$ as in Theorem \ref{th:uniqueness}, $F_K$ is continuous; and if $\cap_j D_j = K$, then $F_{D_j}$ converge uniformly on compact sets of $\CC\PP^n\setminus K_{\rho}$
	to $F_K$. In particular, $F_K$ maps  $\CC\PP^n\setminus K_{\rho}$ onto $\CC\PP^n \setminus K$.
\end{corollary}

\begin{proof} Let $f_c(\zeta) = \rho(c, K) \{(1, c_2,\ldots, c_n) \zeta^{-1} + (1, \bar{c}_2,\ldots, \bar{c}_n) \zeta\} + a_0(c, K)$ be the affine 
representation of an extremal curve as in (\ref{eqn:fc}). We first claim that the function $a_0:  H_{\infty} \ni c \rightarrow a_0(c, K) \in \RR^n$ is continuous.
Suppose to the contrary that there is a sequence $c_j$ such that $c_j \rightarrow c \in H_{\infty}$ and such that $a_0(c_j, K)$ does not converge to $a_0(c, K)$.
Note first that the vectors $\{a_0(c_j, K), \ j = 1,\ldots\} \subset K$ so that without loss of generality we will assume our sequence has $a_0(c_j, K) \rightarrow
\tilde{a}_0 \in K, \ \tilde{a}_0 \neq a_0(c, K).$ Recall from \cite{sic} or Theorem \ref{th:robincont} that $\rho(c, K)$ is continuous in $c$, and so $f_{c_j}(\zeta)$ converge uniformly as maps from $\triangle$
to $\CC\PP^n$ to 
$$\tilde{f}_c(\zeta) \equiv   \rho(c, K) \{(1, c_2,\ldots, c_n) \zeta^{-1} + (1, \bar{c}_2,\ldots, \bar{c}_n) \zeta\} + \tilde{a}_0.$$ It is obvious that 
$\tilde{f}_c \neq f_c$, but that $\tilde{f}_c \in  \mathcal{H}_{(c, K)}(\triangle)$, with $\rho(\tilde{f}_c) = \rho(c, K)$. This contradicts the uniqueness of the
extremal $f_c \in  \mathcal{H}_{(c, K)}(\triangle),$ so $a_0(c, K)$ is continuous in $c$, and therefore $F_K(c/\zeta)$ is continuous as well.

Next, note that $f_{(c, D_j)} \rightarrow f_{(c, K)}$ uniformly on compact sets. This follows from a normal family 
argument and the uniqueness of the extremal in $ \mathcal{H}_{(c, K)}(\triangle)$. Indeed, we have
$f_{(c, D_j)}(\zeta) =  a_1(c, D_j) \zeta^{-1} + a_0(c, D_j) + \sum_{k\geq 1}b_k(c, D_j)\zeta^k$ where $a_1(c,D_j)=\rho(c,D_j)(1,c)$. Since  $\rho(c,D_j)\downarrow \rho(c,K)$ and $f_{(c, D_j)}(\partial \triangle)\subset \bar {D}_j$, the functions $\{a_0(c, D_j)+\sum_{k\geq 1}b_k(c,
D_j)\zeta^k\}_{j=1,...}$ are holomorphic on
$\triangle$, and continuous and uniformly bounded on $\bar \triangle$. Hence $\{f_{(c, D_j)}\}$ form a normal family and any normal limit $\tilde f $
belongs to $\mathcal{H}_{(c, K)}(\triangle)$. But then 
$\rho(c,D_j)\downarrow \rho(c,K)$ and the uniqueness of $f_{(c, K)}$ imply that $\tilde f = f_{(c, K)}$. 
From this and the formula
$$a_0(c, D_j) = \frac{1}{2\pi i} \int_{|\zeta| = \frac{1}{2}} \frac{f_{(c, D_j)}(\zeta)}{\zeta} \, d\zeta,$$
it follows that $a_0(c, D_j) \rightarrow a_0(c, K)$. More precisely, however, we can say that the functions $a_0(c, D_j)$ converge uniformly to $a_0(c, K)$. If 
not, there is an $\epsilon > 0$, and a sequence $c_j \in H_{\infty}$ such that $$|a_0(c_j, D_j) - a_0(c_j, K)| > \epsilon.$$ We can assume $c_j \rightarrow c \in
H_{\infty}$, and that $a_0(c_j, D_j) \rightarrow \tilde{a}_0 \neq a_0(c, K)$, where we have used the continuity of $a_0(c, K)$ as a function of $c$. Then consider
$f_j:\triangle \rightarrow \CC\PP^n$ given in affine coordinates by 
$$f_j(\zeta) = f_{(c_j, D_j)}(\zeta) = a_1(c_j, D_j) \zeta^{-1} + a_0(c_j, D_j) + \sum_{k\geq 1}b_k(c_j, D_j)\zeta^k$$
$$= a_1(c_j, D_j) \zeta^{-1} + a_0(c_j, D_j) +
\bar{a_1(c_j, D_j)}
\zeta + g_j(\zeta),$$ where $g_j(\zeta)$ is holomorphic and uniformly bounded on $\triangle$, independent of $j$, and $g_j(0) = 0$, for all $j$. By what we have
already shown, or assumed,  $a_1(c_j, D_j) \rightarrow a_1(c, K) = \rho(c, K)(1, c_2,\ldots, c_n),$ and $a_0(c_j, D_j) \rightarrow \tilde{a}_0\in K$, and
the
$g_j(\zeta)$ form a normal family on $\triangle$. We will assume that $g_j$ converge uniformly on compact sets in $\triangle$ to a bounded holmorphic function
$\tilde{g}(\zeta)$. Consider the equation $$f_j(e^{i\theta}) = 2\Re(a_1(c_j, D_j) e^{-i\theta}) + a_0(c_j, D_j) + g_j(e^{i\theta}).$$ Taking imaginary parts, we
get $$\Im(g_j(e^{i\theta}))+\Im (a_0(c_j, D_j)) = \Im(f_j(e^{i\theta})).$$ Since $D_j \rightarrow K \subset \RR^n$, we have that $\Im(f_j(e^{i\theta})) \rightarrow
0$, uniformly in
$\theta$. As a result, the harmonic functions $\Im(g_j(\zeta))$, which converge uniformly to $\Im(\tilde{g}(\zeta))$ on compact subsets of $\triangle$, also
converge uniformly on $\bar{\triangle}$ to $0$. Therefore we conclude that $\tilde{g}$ is a real constant function, and since $0 \equiv g_j(0) \rightarrow
\tilde{g}(0)$, this constant is $0$. Thus the normal limit $\tilde f$ of the $f_j's$ has the form $\tilde{f}(\zeta) = a_1(c, K) \zeta^{-1} + \tilde{a}_0 +\bar{a_1(c, K)} \zeta$. The limit procedure implies
$\tilde{f}(0) = c$, and $\tilde{f} \neq f_{(c, K)}$ since $\tilde{a}_0 \not = a_0(c,K)$. Furthermore, for each fixed $\zeta \in \triangle,$ given any $\delta > 0$, we have $$|V_K(\tilde f(\zeta)) -
\log \frac{1}{|\zeta|}| < |V_K(f_{(c_j, D_j)}(\zeta)) - \log \frac{1}{|\zeta|}| + \delta$$ $$ < |V_{\bar {D}_j}(f_{(c_j, D_j)}(\zeta)) - \log \frac{1}{|\zeta|}| + 2
\delta = 2\delta,$$ for all $j$ sufficiently large, since $V_K$ is continuous, $f_{(c_j, D_j)}(\zeta) \rightarrow \tilde{f}(\zeta)$ and $V_{\bar {D}_j} \rightarrow V_K$
uniformly on $\CC^n$. Since $\delta$ is arbitrary, $V_K(\tilde{f}(\zeta)) = \log \frac{1}{|\zeta|},$ for all $\zeta \in \triangle$; since $V_K$ is continuous on
$\CC^n$, and $\tilde{f}$ is continuous on $\bar{\triangle}$, we have $V_K(\tilde{f}(e^{i\theta})) \equiv 0, \ \theta \in [0, 2\pi]$. Thus $\tilde{f}(e^{i\theta})
\in K$, for all $\theta$, showing that $\tilde{f} \in  \mathcal{H}_{(c, K)}(\triangle)$. We conclude that $\tilde f$ is an extremal, contradicting  the uniqueness of the extremal $f_{(c, K)} \in 
\mathcal{H}_{(c, K)}(\triangle)$.

Finally, if $F_{D_j}$ do not converge uniformly to $F_K$, we would have an $\epsilon > 0$ and sequences $c_j \in H_{\infty}, \ \zeta_j \in \triangle$, such that 
$$\mbox{dist}_{\CC\PP^n}(F_{D_j}(c_j/\zeta_j), F_K(c_j/\zeta_j )) > \epsilon, \ j = 1,\ldots,$$ and where $\zeta_j \rightarrow \zeta \in \triangle,
\ c_j \rightarrow c \in H_{\infty}$. But for $j >> 0$, we have as before 
$$f_{(c_j, D_j)}(\zeta_j) = a_1(c_j, D_j) \zeta_j^{-1} + a_0(c_j, D_j) + \bar{a_1(c_j, D_j)} \zeta_j + g_j(\zeta_j),$$ where the $g_j$ are a normal family, and 
so must converge normally to a limit $\tilde{g}.$ By the assumption, we must have $ \lim_{j\to \infty} g_j(\zeta_j) = \tilde{g}(\zeta) \neq 0$, whereas the same argument
concerning $\Im \, g_j$ converging uniformly to $0$ still applies, leading to the conclusion $\tilde{g} \equiv 0$, a contradiction.

\end{proof}

\noindent {\bf Remark.} We have used the basic hypothesis that the extremal curve $f_c \in  \mathcal{H}_{(c, K)}(\triangle)$ is unique. We point out that 
all the arguments above go through for $K \subset \RR^n$ if $K$ is symmetric, i.e., $K = -K$. Although the extremal disk is not necessarily unique in such a case, we may use the
extremal disk $f \in  \mathcal{H}_{(c, K)}(\triangle)$ for each $c \in H_{\infty}$ such that $a_0(f) = 0 \in \RR^n$. That such an extremal disk exists for every $c
\in H_{\infty}$ follows from taking an approximating sequence $D_j$ such that $D_j$ is symmetric, as well as conjugation invariant ($\sigma(D_j) = D_j$), followed by a
normal family argument as above. This implies that $a_0(c, D_j) = 0$, for every $c \in H_{\infty}$. It is clear that such a symmetric extremal is unique, and given
by $f_c(\zeta) = \rho(c, K) \{(1, c_2,\ldots, c_n) \zeta^{-1} + (1, \bar{c}_2,\ldots, \bar{c}_n) \zeta\}.$ However, if the approximating sets are not symmetric ($D_j\not = -D_j$), a limit disk might not satisfy $a_0(f)=0$; indeed, the limit of the sequence of points $a_0(c, D_j)$ need not exist:

\noindent {\bf Example.} Take the square $K=[-1,1]\times [-1,1]$ in $\RR^2$. Fix $0<a<1$, then given a decreasing sequence of positive numbers 
$\epsilon_j\downarrow 0$, define a decreasing sequence of non-symmetric heptagons $K_j$ by adding the vertices $(1+\epsilon_j,a),(-(1+\epsilon_j),a)$ and $(0,1+\epsilon_j)$ to the vertices $(\pm 1,\pm1)$ of the square. We can choose a decreasing sequence $D_j$ of bounded, smoothly bounded
strictly convex open sets in $\CC^n$ with $\sigma(D_j)=D_j$ such that $\cap D_j = K$ with $D_j$ sufficiently close to $K_j$ so that for each $c\in H_{\infty}$ the
limit points of the sequences $a_0(c, D_j)$ and $a_0(c, K_j)$ are the same. For $c=[0:1:0]$ in homogeneous coordinates, the extremal map for $K_j$ is clearly
$f_{(c,K_j)}(\zeta)=(0,a)+(\frac{1+\epsilon_j}{2}(\zeta+1/\zeta),0)$; in particular, $a_0(c,K_j)=(0,a)$. Hence $f_c(\zeta)=(0,a)+(\frac{1}{2}(\zeta +1/\zeta),0)$
and $a_0(f_c)=(0,a)$. Moreover, if for $j$ even we take $K_j$ as described and for $j$ odd we replace the vertices $(1+\epsilon_j,a),(-(1+\epsilon_j),a)$ by
$(1+\epsilon_j,-a),(-(1+\epsilon_j),-a)$, we obtain, provided $\epsilon_j \to 0$ appropriately, an alternating nested sequence of up- and down-going ``coffins''. In
this case, 
$a_0(c,K_j)$ for $c=[0:1:0]$ alternates between $(0,a)$ and $(0,-a)$, giving two limit points.

In section \ref{sec:FhomeoIII}, we discuss a ``canonical'' choice of center function $c\to a_0(c, K)$ for nonsymmetric $K$. 

\vskip5mm
 
\section{\bf $F_K$ is a homeomorphism, case I.}
\label{sec:FhomeoI}

\vskip3mm

In this section we show that in favorable circumstances the extremal curves give a continuous foliation of $\CC\PP^n \setminus K$.

\begin{theorem} 
\label{th:fol1} Let $K\subset \subset \RR^n$ be convex such that the extremal curve through any $c \in H_{\infty}$ is unique. Then the Robin exponential map $F_K$ is a   homeomorphism of $\CC\PP^n \setminus {K_{\rho}}$ onto $\CC\PP^n \setminus K$.
\end{theorem}

\begin{proof} We have already shown that $F_K$ is continuous, and since $$F_K^{*}V_K (c/\zeta)= V_K(f_c(\zeta))=\log \frac{1}{|\zeta|}$$ on $\CC^n \setminus K$, 
we  have that $F_K$ is a proper mapping. It suffices, therefore, to show that $F_K$ is one-to-one. Since each $f_c$ is an embedding of $\triangle$ into $\CC\PP^n
\setminus K$, we have to show that there do not exist two distinct points $c, \tilde{c} \in H_{\infty}$, and values $\zeta, \tilde{\zeta} \in \triangle$ such that $z_0 =
f_c(\zeta) = f_{\tilde{c}}(\tilde{\zeta})$. Since $f_c(0) = c \neq \tilde{c} = f_{\tilde{c}}(0),$ we see that both $\zeta, \tilde{\zeta} \neq 0,$ and hence $z_0
\in \CC^n$. Let us consider in place of the disk $f_c(\triangle)$ the projective curve $C := f_c(\CC\PP^1)$ extending it, and similarly for $\tilde{C} :=
f_{\tilde{c}}(\CC\PP^1)$. Such a  curve is a projective line (doubly covered) or a non-singular quadric depending on whether the point $c$ (or $\tilde{c}$) is real
or not, respectively. We will show that:
$$\mbox{for extremal curves} \, C, \tilde{C}, \ C \cap \tilde{C} \subset K.$$ The proof of this, which will prove the theorem, is an exercise in elementary
geometry,  but seems to require considering several cases, depending on whether $C$ and $\tilde{C}$ are linear or quadratic, and the dimension of the linear span
of $C \cup
\tilde{C} \subset \CC\PP^n$. First note that since $\sigma(C) = C$, and similarly for $\tilde{C}$, the set $C \cap \tilde{C}$ is sent to itself by $\sigma$. Let
$C_{\RR}$ denote the real points of $C$, and similarly for $\tilde{C}$.

\vskip 2mm

\noindent \emph{Case: n = 2.} 

\vskip 2mm
Within this case we will treat three subcases: a) $C, \tilde{C}$ both linear, b) $C$ linear and $\tilde{C}$ quadric, and c) both $C, \tilde{C}$ are quadrics.

\vskip 2mm

\noindent a) In this case, $C \cap \tilde{C} = \{z_0\}$ is one point, which must be real since the set $C \cap \tilde{C}$ is $\sigma$-invariant. 
If $z_0$ is not in $K$, then $C_{\RR} \cap K$ is a maximal length segment within $K$ in its direction, and similarly for $\tilde{C}_{\RR}$, and these two segments
do not cross within $K$. If we denote the two endpoints of $C_{\RR} \cap K$ by $a, b$, and the endpoints of $\tilde{C}_{\RR} \cap K$ by $c, d$, we can assume that
$b$ is between $a$ and $z_0$, and similarly $c$ is between $d$ and $z_0$, then the convex hull of $(C_{\RR} \cap K )\cup (\tilde{C}_{\RR} \cap K)$ is the
quadrilateral
$\mathcal{H}$ bound by $\bar{ab}, \bar{bc}, \bar{cd}$, and $\bar{da}$, and is contained entirely within $K$. If the two segments $\bar{bc}$ and $\bar{da}$ are not
parallel segments in $\RR^2$, then one or the other of the seg\-ments $\bar{ab}$ or $\bar{cd}$ can be deformed parallelly within $\mathcal{H}$ in such a way as to
in\-crease its length, contradicting the extremality of $C$ or $\tilde{C}$. If the segments $\bar{bc}$ and $\bar{da}$ are parallel, then both of the segments
$\bar{ab}$ and $\bar{cd}$ can be deformed parallelly within the convex hull $\mathcal{H}$  in such a way as to preserve their lengths, contradicting the uniqueness
of the curves $C, \tilde{C}$ as extremals. Thus, we conclude that $z_0 \in K$.

\vskip 2mm

\noindent b) In this case, $C \cap \tilde{C}$ is two points, counted with multiplicity. Let us assume that the segment $\bar{ab} = C_{\RR} \cap K$ is disjoint 
from $\tilde{C}_{\RR} \subset K$. (If they intersect, they obviously intersect in two points, counted with multiplicities.) Again, we consider the convex hull
$\mathcal{H}$ of $\bar{ab} \cup \tilde{C}_{\RR}$. There are two points $c, d \in \tilde{C}_{\RR}$ such that $\mathcal{H}$ is bound by the line segments $\bar{ab},
\bar{bc}, \bar{da}$ and an arc $\bar{\bar{cd}} \subset \tilde{C}_{\RR}.$ If the segments $\bar{bc}, \bar{da}$ are not parallel, then as above, one of $\bar{ab},
\tilde{C}$ cannot be of maximal length or area in its class of segments or ellipses. If the segments are parallel, then $C$ cannot be a unique extremal.
Thus we conclude that $\bar{ab} \cap \tilde{C}_{\RR} \neq \phi$, and therefore, $C \cap \tilde{C} \subset K$. 

\vskip 2mm

\noindent c) In this case, we have two ellipses $C_{\RR}, \tilde{C}_{\RR} \subset K$ and we have to show that $C_{\RR} \cap \tilde{C}_{\RR}$ consists of four 
points, counting multiplicities. The possibilities are $0, 2$ or $4$. We cannot have one of the ellipses strictly contained within the other, since then the inner
one could not be of maximal area in its family. Therefore, if $C_{\RR} \cap \tilde{C}_{\RR}$ is 0 or 2 points counting multiplicities, then there must be four
points $a, b \in C_{\RR}, \ c, d \in \tilde{C}_{\RR}$ so that the convex hull $\mathcal{H}$ of $C_{\RR} \cup \tilde{C}_{\RR}$ must be bound by an arc
$\bar{\bar{ab}}
\subset C_{\RR},$ an arc $\bar{\bar{cd}} \subset \tilde{C}_{\RR}$, and two segments $\bar{bc}, \bar{da}$. If the segments $\bar{bc}, \bar{da}$ are not parallel,
then once again, one of $C_{\RR}, \tilde{C}_{\RR}$ cannot be extremal in its family. If the segments $\bar{bc}, \bar{da}$ are parallel, neither $C_{\RR}$ nor
$\tilde{C}_{\RR}$ can be a unique extremal in its family. 

\vskip 2mm

\noindent This concludes the proof in the case $n = 2$.

\vskip 2mm
\noindent \emph{Case: $n > 2$.}

\vskip 2mm

Within this case, the number of points of intersection in $C \cap \tilde{C}$ is not determined beforehand.
We will treat subcases a), b) and c) again as above. 

\vskip 2mm

\noindent a) If $C, \tilde{C}$ are both linear, they intersect in 1 point or none. If none, we are done. If one, then the linear hull of 
$C \cup \tilde{C} \subset \CC\PP^n$ is a projective plane $\Pi$, and $\sigma(\Pi) = \Pi$. Since the segments $C_{\RR} \cap K, \ \tilde{C}_{\RR} \cap K$ are extremal
for their directions in $K$, they will also be extremal for the convex set $K \cap \Pi := K_{\Pi}$, and we reduce the argument to subcase a) above in Case $n = 2$.

\vskip 2mm

\noindent b) In this case, $C \cap \tilde{C}$ can consist of $0, 1,$ or $2$ points, counting multiplicity. If 0, we are done, and if one, then this point of 
intersection must be in $\RR^n$, and so in $\tilde{C}_{\RR} \subset K$, which is what we want to prove. Finally, if there are two points of intersection, then $C,
\tilde{C}$ are coplanar, contained in a plane $\Pi$. Since the segment $C_{\RR} \cap K$ and the ellipse $\tilde{C}_{\RR}$ are extremal for $K$, they are also for
$K \cap \Pi \subset \Pi$. This reduces the question to subcase b) of Case $n =2$.

\vskip 2mm

\noindent c) In this case, the possibilities for $C \cap \tilde{C}$ are 0, 1, 2, 3, or 4 points (with multiplicities). If we have 0, we are done, and if we have 
1, then it must be a real point, and so in $K$. If we have 3 or 4, then $C$ and $\tilde{C}$ are coplanar, contained in a plane $\Pi$. Again, $C_{\RR}$ and
$\tilde{C}_{\RR}$ are extremal for $K \cap \Pi$ and we conclude $C \cap \tilde{C} \subset K \cap \Pi \subset K$, by subcase c) of Case $n = 2$. If we have 2 points
of intersection, they must be either 2 real points, and we are done, or two conjugate, non-real points. We must show that this latter case cannot occur, for
$C_{\RR}, \ \tilde{C}_{\RR}$ extremal in their families in $K$. Note that the projective hull of $C$ is a plane $\Pi$, and similarly for $\tilde{C} \subset
\tilde{\Pi}$, and
$\Pi \cap \tilde{\Pi}$ is the projective line determined by the two (unequal) non-real points of intersection in $C \cap \tilde{C}$. This line is real, i.e.,
$\sigma$-invariant.

This case will require some explicit computation, and we use real affine diffeomorphisms of $\RR^n$ (and therefore of $\CC^n, \CC\PP^n$) to simplify the situation. 
Note that this is possible, since such a real affine transformation ${\mathcal A}$ will take $K$ to another real convex body in $\RR^n$, and will take lines and quadrics
to lines and quadrics. Furthermore, it is easy to check that it also takes a curve $C$ extremal for $K$ in the family determined by $c \in H_{\infty}$ to ${\mathcal A}(C)$
which will be extremal for ${\mathcal A}(K)$ in the family determined by ${\mathcal A}(c) \in H_{\infty}$. Thus, we can assume that we have real coordinates $x = (x_1,\ldots, x_n)$ such that $C_{\RR}$ is an ellipse contained in the plane of $x_1, x_2$, while $\tilde{C}_{\RR}$ is contained in the plane of $x_1, x_3$, and the quadrics $C,
\tilde{C}$ intersect in two non-real points $(t, 0,\ldots, 0) \neq (\bar t, 0,\ldots, 0)$; i.e., $t =u+iv$ with $v\not = 0$.

We may also scale the coordinates by a real linear transformation ${\mathcal A}: (x_1, x_2, x_3, x_4,\ldots, x_n) \rightarrow (x_1, \lambda_2 \cdot x_2,
\lambda_3 \cdot x_3, x_4,\ldots, x_n),$ where $\lambda_2, \lambda_3 > 0$. If $\lambda_2, \lambda_3$ are suitably chosen, ${\mathcal A}(C_{\RR})$ and ${\mathcal A}(\tilde{C}_{\RR})$ will
be circles, i.e., ellipses of eccentricity $e = 1$, so we will henceforth assume they are both circles. We can also consider the problem in $\CC^3 = \{(z_1, z_2,
z_3, 0,\ldots, 0)\} \subset \CC^n$, since $C_{\RR}, \tilde{C}_{\RR}$ will be extremal also for $K \cap \CC^3$.

So, we have two circles

$$C_{\RR} = \left\{\begin{array}{l}
  (x_1 - \alpha)^2 + (x_2 - \beta)^2 = r^2, \ \alpha, \beta, r \in \RR, \ r > 0,\, \mbox{and} \\
  x_3 = 0.
  \end{array} \right. $$
and

$$\tilde{C}_{\RR} = \left\{\begin{array}{l}
  (x_1 - \gamma)^2 + (x_3 - \delta)^2 = \tilde{r}^2, \ \gamma, \delta, \tilde{r} \in \RR, \  \tilde{r} > 0,\, \mbox{and} \\
  x_2 = 0.
  \end{array} \right. $$

Now if $C \cap \tilde{C}$ is to contain two solutions, we must have that the two complex equations

$$(z_1 - \alpha)^2 + \beta^2 = r^2 \ \hbox{and} \ (z_1 - \gamma)^2 + \delta^2 = \tilde{r}^2$$
are the {\em same} equation.  This is if and only if we have 
$$\alpha = \gamma \ \hbox{and} \ \beta^2 - r^2 = \delta^2 -\tilde{r}^2.$$
By renumbering coordinates if necessary, we can assume $\delta \geq \beta$, and by translating in the $x_1$ direction, we can assume $\alpha = \gamma = 0$. Thus
$$t ^2 + \beta^2 = r^2 \ \hbox{and} \ t^2 + \delta^2 = \tilde{r}^2.$$
The condition that
$$\bar t^2 + \beta^2 = r^2 \ \hbox{and} \ \bar t^2 + \delta^2 = \tilde{r}^2$$
as well implies that $u=0$; then $(\pm iv,0,...,0)\in C\cap \tilde C$ implies that
$$-v^2 =r^2-\beta^2=\tilde r^2 -\delta^2 <0$$
so that $\beta^2 >r^2$ and $\delta^2 > \tilde r^2$. Thus we may assume, for simplicity, that $x_2 > 0$ on $C_{\RR}$ and $x_3 > 0$ on
$\tilde{C}_{\RR}$. 

Let us first treat the case $\delta = \beta$, and therefore $r = \tilde{r}$. In this case, we can calculate directly that the convex hull 
$\mathcal{H}(\beta) \subset K$ of $C_{\RR} \cup \tilde{C}_{\RR}$ is the intersection of the sector $S = \{x_2 \geq 0, x_3 \geq 0\}$ with the cylindrical region
bounded by the image of the map $$B: [0, 2\pi] \times [0,1] \ni (\theta, t) \rightarrow (r \cos \theta, t(\beta + r\sin \theta), (1-t)(\beta +
r\sin
\theta)) \in \RR^3.$$ But this implies that for $\epsilon > 0$ sufficiently small, the plane $\{x_3 = \epsilon\} \cap \mathcal{H}(\beta)$ is a circle
$(r \cos \theta, \beta + r\sin \theta-\epsilon, \epsilon) $ congruent to $C_{\RR}$ and in its family, contradicting that $C_{\RR}$ is the unique extremal for $K$
in its family. 

Finally, in the case that $\delta > \beta$, we have that $\delta - \tilde{r} < \beta - r$ and that $\delta + \tilde{r} > \beta + r.$ This in turn implies 
the convex hull of $C_{\RR} \cup \tilde{C}_{\RR}$ contains the set $\mathcal{H}(\beta)$. This shows that $C_{\RR}$ is not extremal within its family, a
contradiction. We thus conclude that $C$ and $\tilde{C}$ cannot intersect in two non-real points, completing the proof of Theorem \ref{th:fol1}.

\end{proof}

\noindent {\bf Remark}. The proof actually shows more: if {\it one} of $C$ or $\tilde{C}$ is a unique extremal for its family, then $C\cap \tilde C\subset K$. This observation will be used in section 10.

\vskip5mm

 \section{\bf $F_K$ is a homeomorphism, case II.}

 \label{sec:FhomeoII}

\vskip 3mm
 
We say that a convex body $K\subset \RR^n$ has property ${\mathcal Q}$, and we write $K\in {\mathcal Q}$, if $K$ has unique extremals; i.e., $K$ satisfies 
the conclusion of Theorem \ref{th:uniqueness}. It is clear that given any $K\not \in {\mathcal Q}$ we can approximate $K$ from outside by a sequence $\{K_j\}\subset {\mathcal Q}$
with $K_{j+1}\subset K_j$ and $\cap K_j =K$ (we write $K_j \downarrow K$). We show that this can be done in such a way that, as in Theorem 7.1, the mappings
$F_{K_j}(c/\zeta)=a_0(c,K_j)+\rho(c,K_j)\{(1,c_2,...,c_n)\zeta^{-1} +(1,\bar c_2,...,\bar c_n)\zeta\}$ (see (\ref{eqn:fkc})) converge locally uniformly to a map $F_K$ from $\CC\PP^n \setminus K_{\rho}$ (recall (\ref{krho})) onto $\CC\PP^n \setminus K$. We observe that $K_j \downarrow K$ implies that $V_{K_j}\to  V_K$ uniformly on $\CC^n$ and $\rho(c,K_j)\to
\rho(c,K)$ uniformly on $H_{\infty}$.

We begin with a preliminary result.

\begin{lemma}
	\label{l:1}
	Let $K\subset \subset \RR^n$ be a convex body and suppose $\{K_j\}\subset {\mathcal Q}$ with $K_j \downarrow K$. If the limit
$$a(c):=\lim_{j\to \infty}a_0(c,K_j)$$
exists, and the limit is uniform on $H_{\infty}$, then
$$F_K(c/\zeta)=a(c)+ \rho(c,K)\{(1,c_2,...,c_n)\zeta^{-1} +(1,\bar c_2,...,\bar c_n)\zeta\}$$
maps $\CC\PP^n \setminus K_{\rho}$ onto $\CC\PP^n \setminus K$.
\end{lemma}

\begin{proof} The proof follows along lines similar to Corollary \ref{cor:uniqueness}. We first show that for each $c$,
$$f(\zeta):=a(c)+\rho(c,K)\{(1,c_2,...,c_n)\zeta^{-1} +(1,\bar c_2,...,\bar c_n)\zeta\}$$
is an extremal map for $K$. By hypothesis and the fact that $\rho(c,K_j)\to
\rho(c,K)$, if $c^{(j)}\to c$,
$$f_{(c^{(j)},K_j)}(\zeta)=a_0(c^{(j)},K_j)+\rho(c^{(j)},K_j)\{(1,c^{(j)}_2,...,c^{(j)}_n)\zeta^{-1} +(1,\bar c^{(j)}_2,...,\bar c^{(j)}_n)\zeta\}$$ converges uniformly on compact subsets of $\triangle$ to $f(\zeta)$. 
Since $V_{K_j}\to V_K$ uniformly on $\CC^n$, given $\delta >0$, for each fixed $\zeta \in \triangle$, if $j\geq j_0(\delta)$ is sufficiently large, 
$$|V_K( f(\zeta)) -
\log \frac{1}{|\zeta|}| < |V_K(f_{(c^{(j)}, K_j)}(\zeta)) - \log \frac{1}{|\zeta|}| + \delta$$ $$ < |V_{K_j}(f_{(c^{(j)}, K_j)}(\zeta)) - \log \frac{1}{|\zeta|}| + 2
\delta = 2\delta.$$
Thus $V_K( f(\zeta))=\log \frac{1}{|\zeta|}$ on $\triangle$; by continuity, $V_K(f(e^{i\theta}))=0$ and $f \in \mathcal{H}_{(c, K)}(\triangle)$. Hence $f$ is
an extremal for $K,c$.

Since $\rho(c,K_j)\to
\rho(c,K)$ uniformly on $H_{\infty}$, uniform convergence of the center functions $a_0(c,K_j)$ to $a(c)$ is equivalent to local uniform convergence of
$F_{K_j}(c/\zeta)$ to
$F_{K}(c/\zeta)$, yielding the conclusion.
\end{proof}

The argument above yields uniform convergence of $f_{(c^{(j)},K_j)}$ to $f$ on $\bar \triangle$; a fact we will utilize repeatedly. We therefore record the 
precise statement needed as a corollary.

\begin{corollary}
	\label{c:1}
	Let $K\subset \subset \RR^n$ and $\{K_j\}\subset {\mathcal Q}$ with $K_j \downarrow K$ satisfying the hypothesis of the lemma. If for a subsequence of positive 
integers $\{j_k\}$, we have $a:=\lim_{k\to \infty}a_0(c^{(j_k)},K_{j_k})$ exists and  $c:=\lim_{k\to \infty}c^{(j_k)}$ exists, then
$$f(\zeta):=a+\rho(c,K)\{(1,c_2,...,c_n)\zeta^{-1} +(1,\bar c_2,...,\bar c_n)\zeta\}$$
is an extremal in $\mathcal{H}_{(c, K)}(\triangle)$.
\end{corollary}

We next verify that $F_K$ as described above is a homeomorphism if $n=2$.

\begin{proposition}
	\label{prop:1}
	Let $K\subset \subset \RR^2$ and $\{K_j\}\subset {\mathcal Q}$ with $K_j \downarrow K$ satisfy the hypothesis of the lemma. Then $F_K$ is a homeomorphism. 
\end{proposition}

\begin{proof} As in the proof of Theorem 8.1, we must show that if $c,\tilde c$ are distinct points in $H_{\infty}$, then the extremal curves $C:=f_c(\CC\PP^1)$
and $\tilde C:=f_{\tilde c}(\CC\PP^1)$ satisfy $C\cap \tilde C\subset K$. We use an intersection multiplicity argument. Suppose there exists a point $p\not\in
K$ with $p\in C\cap \tilde C$. We claim that this must be a transverse intersection. For if not, since $V_K(p)=\log r$ for some $r>1$, by rotating coordinates in $\triangle$ if necessary, we can find $\zeta_0\in \triangle$ with
$f_c(\zeta_0)=f_{\tilde c}(\zeta_0)=p$ and $f_c'(\zeta_0)=f_{\tilde c}'(\zeta_0)$. Let $h(\zeta):=f_c(\zeta)-f_{\tilde c}(\zeta)$. Then $h$ is of the form 
$$h(\zeta)=\alpha +\beta /\zeta +\bar \beta \zeta$$
with $\alpha \in \RR^2$, $\beta\in \CC^2$, and $h(\zeta_0)=h'(\zeta_0)=0$. But $h'(\zeta_0)=0$ gives $|\zeta_0|=1$ which is impossible since $p\not\in K$.

Fix a small ball $B$ containing $p$ which is disjoint from $K$. Then for large $j$, $B\cap K_j=\emptyset$. Since $f_{(c, K_j)},f_{(\tilde c, K_j)}$ converge 
uniformly to $f_c,f_{\tilde c}$ in a neighborhood of $\zeta_0$, $C_j:=f_{(c, K_j)}(\CC\PP^1)$ must intersect $\tilde C_j:=f_{(\tilde c, K_j)}(\CC\PP^1)$ in $B$ for $j$
large (moreover, transversally), a contradiction.
\end{proof}

We conclude this section with a construction in $\RR^2$ satisfying the hypothesis of Lemma 9.1. Let $(x,y)$ be coordinates in $\RR^2$. For the remainder of this section, by abuse of notation, we write $c=[0:1:c]\in H_{\infty}$. 

\begin{proposition}
	\label{prop:2}
	Let $K\subset \subset \RR^2$ be a convex body with $K\not \in {\mathcal Q}$. Then there exist $\{K_j\}\subset {\mathcal Q}$ with $K_j \downarrow K$ such that 
the limit 
$$a(c):=\lim_{j\to \infty}a_0(c,K_j)$$
exists, and the limit is uniform on $H_{\infty}$. In particular, 
$$F_K(c/\zeta)=a(c)+\rho(c,K)\{(1,c)\zeta^{-1} +(1,\bar c)\zeta\}$$
is a homeomorphism of $\CC\PP^2 \setminus K_{\rho}$ onto $\CC\PP^2 \setminus K$.
\end{proposition}

\begin{proof} Let
$$G:=\{c\in H_{\infty}: \hbox{the extremal curve for} \ K \ \hbox{through} \ c \ \hbox{is unique}\}.$$
We first show that for {\it any} approximation $K_j \downarrow K$ with  $\{K_j\}\subset {\mathcal Q}$, 
$$a(c):=\lim_{j\to \infty}a_0(c,K_j)$$
exists for $c\in G$ and the convergence is uniform; i.e., if $\{c_j\}\in H_{\infty}$ with $c_j\to c$, then $a_0(c_j,K_j)\to a(c)$. First we verify existence of 
the limit: take any subsequential limit, say $\tilde a$, of the full sequence $\{a_0(c_j,K_j)\}$. Consider
$$\tilde f(\zeta):=\tilde a+ \rho(c,K)\{(1,c)\zeta^{-1} +(1,\bar c)\zeta\}.$$
By the argument in Lemma 9.1, i.e., by Corollary 9.1, since $K_j \downarrow K$ implies $V_{K_j}\to V_K$ uniformly and $\rho(\cdot,K_j)\to \rho(\cdot,K)$ uniformly, we
see that
$\tilde f$ is an extremal for $K,c$. Since $c\in G$, $\tilde f=f_c$ is unique; i.e., there is a unique limit point $a(c)$ of the sequence $\{a_0(c_j,K_j)\}$. 

Now $K$ can have at most countably many pairs $\{L_j,M_j\}$ of parallel segments in $\partial K$. Fix one such pair $ L_j,M_j$; we write $L:=L_j, \ M:=M_j$. By rotating and translating 
coordinates, we may assume that $L$ is the horizontal segment $[0,2p], \ p>0$, on the positive $x-$axis, and $M$ is contained in a horizontal line
$y=y_0<0$. Under the assumption that the $x-$axis is not a horizontal tangent to $\partial K$ at $(0,0)$ or $(2p,0)$, we can deform $L$ via a linear midpoint
homotopy: define, for $0\leq \alpha \leq 1$,
$$H(\alpha,t)= t(2p,2\alpha p), \ \hbox{for} \ 0\leq t\leq 1/2;$$
$$H(\alpha,t)=2p(t,\alpha(1-t)),
 \ \hbox{for} \ 1/2\leq t\leq 1.$$  
We perform this construction on each side $L_j, \ j=1,2,...$; we call the corresponding homotopty $H_j$. Now choose $\alpha^1_1>0$ sufficiently small so that replacing $L_1$ by 
$\tilde L_1=\{H_1(\alpha_1^1,t):t\in [0,1]\}$ in $\partial K$, we introduce no new parallel segment and still have a convex set; then choose $\alpha^1_2>0$
sufficiently small so that replacing $L_2$ by $\tilde L_2=\{H_2(\alpha^1_2,t):t\in [0,1]\}$ in $\partial K$, we introduce no new parallel segment and still have a
convex set; continue choosing $\alpha^1_j>0$ sufficiently small so that replacing $L_j$ by $\tilde L_j=\{H_j(\alpha^1_j,t):t\in [0,1]\}$ in $\partial K$, we
introduce no new parallel segment and still have a convex set, $j=1,2,...$; this yields the set $K_1$. Next we choose $\alpha^2_j <\alpha^1_j$, $j=1,2,...$
successively to construct $K_2$; and we continue this process to obtain a sequence $\{K_j\}\subset {\mathcal Q}$ with $K_j \downarrow K$. 
\medskip

\noindent {\bf Remark}. The use of a linear
homotopy, let alone of the ``midpoint'' version, is not at all essential; a continuous homotopy with the ``tops'' converging to a unique point on $L$ will suffice
provided the resulting $K_j$ are convex and satisfy ${\mathcal Q}$. 
\medskip

Note that $H_{\infty}\setminus G=\cup_{j=1}^{\infty} B_j$ where
$$B_j:=\{c\in H_{\infty}\setminus G: \hbox{extremal} \ C \ \hbox{for} \ c \ \hbox{has} \ C\cap K \ \hbox{hits both} \ L_j,  \ M_j\}.$$
Moreover, this is a disjoint union of open sets; i.e., $B_j\cap B_k=\emptyset$ for $j\not = k$. For each $c\in H_{\infty}\setminus G$, denote the segment of 
centers of extremals for $K,c$ by $[a_1(c),a_2(c)]$ and the ``tops'' of these extremals on the corresponding segment $L$ by $[T_1(c),T_2(c)]$. In the arguments
below, we tacitly assume we are working on a subset $B:=B_j$ of $H_{\infty}$ with coordinates chosen so that $L:=L_j, \ M:=M_j$ are horizontal; and $L=[0,2p]$. For use in the next section, we prove the following.

\begin{lemma}
	\label{topcont}
The map $c\to a_1(c)$ $(c\to a_2(c))$ is continuous on $H_{\infty}\setminus G$. Hence, the map $c\to T_1(c)$ $(c\to T_2(c))$ is continuous.
\end{lemma}

\begin{proof} Fix $a\in (a_1(c),a_2(c))$. Then for $c'\in B$ sufficiently near $c$, the ellipse $E_{c'}:=\{a+ \rho'
\{(1,c')\zeta^{-1} +(1,\bar c')\zeta\}: 
|\zeta|=1\}$, where $\rho'$ is chosen so that $E_{c'}$ has its top $t$ (here and below, ``top'' will mean the point with largest $y-$coordinate) in
$(T_1(c),T_2(c))$ and its bottom $b$ (point with smallest $y-$coordinate) on the segment $M$, is contained in $K$; we show it is an extremal for $K,c'$; i.e., $\rho' =\rho(c',K)$. This follows
by convexity: if not, $\rho'<\rho(c',K)$ and for any extremal ellipse for $K,c'$ with top $t'$ and bottom $b'$, the segment joining $t$ to $t'$ and
the segment joining $b$ to $b'$ lie in $K$. Since $\rho'<\rho(c',K)$, at least one of these segments is not horizontal, contradicting the assumption that $L,M$
are parallel segments in $\partial K$.

Hence for any $a\in (a_1(c),a_2(c))$, $a> a_1(c')$ for $c'\in B$ sufficiently close to $c$. Hence $\limsup_{c'\to c, \ c'\in B}a_1(c')\leq a_1(c)$. Let 
$b(c):=\liminf_{c'\to c, \ c\in B}a_1(c')$. Taking a sequence $\{c_j\}$ with $c_j \to c$ so that $a_1(c_j)\to b(c)$, utilizing Corollary 9.1, we conclude that $f(\zeta):=
b(c)+\rho(c,K)\{(1,c)\zeta^{-1} +(1,\bar c)\zeta\}$ is an extremal for $K,c$ so that $b(c)\in [a_1(c),a_2(c)]$ and hence $\lim_{c'\to c, \ c'\in B}a_1(c')=a_1(c)$.
\end{proof}

To finish the proof of the proposition, we need verify the uniform convergence of the sequence of functions $a_0(c,K_j)$. To do this we must show that for $c\in B$, if $\{c_j\}\in B$ with $c_j \to c$, then $\lim_{j\to \infty} a_0(c_j,K_j)$ exists. 
For each $c\in H_{\infty}$, we write $E_j(c)$ for the ellipse $\{f_{(c, K_j)}(\zeta): |\zeta|=1\}$ in $K_j$. We consider two cases according to whether the midpoint $p$ of $L$ lies within the open ``top interval'' $(T_1(c),T_2(c))$ or not.
\medskip

\noindent {\sl Case I: $p\in (T_1(c),T_2(c))$}: We have $E_j(c_j)$ is the ellipse $f_{(c_j,K_j)}(\partial \triangle)$ where $f_{(c_j,K_j)}(\zeta)= a_0(c_j,K_j)+ 
 \rho(c_j,K_j)\{(1,c_{j})\zeta^{-1} +(1,\bar c_{j})\zeta\}$.  We show that 
$$\lim_{j\to \infty} a_0(c_j,K_j)=:a$$ exists; moreover, we show that the ellipse 
$$\{a+  \rho(c,K)\{(1,c)\zeta^{-1} +(1,\bar c)\zeta\}: |\zeta|=1\}$$
has top $t=p$. Note that for $j$ large, $E_j(c_j)$ must hit $L$ (and generally in two points); let $l_j$ denote such a point if it lies to the left of $p$ 
($x-$coordinate of $l_j$ less than $p$) and let $r_j$ denote such a point if it lies to the right of $p$. We take $j$ large so that at least one of $l_j,r_j$
exists. Suppose, for the sake of obtaining a contradiction, that a subsequence of the $l_j'$s converges to $\tilde p <p$. Form the ellipse $E$ with top at $\tilde
p$ and bottom $\tilde b$ on $M$ using the scale factor $\rho(c,K)$; i.e., 
$$E=\{f(\zeta):=\tilde a+ \rho(c,K)\{(1,c)\zeta^{-1} +(1,\bar c)\zeta\}: |\zeta|=1\}$$
where $\tilde a$ is the midpoint of $\tilde p$ and $\tilde b$. Corollary 9.1 implies that $E$ is an extemal ellipse for $K,c$; moreover, by the uniform convergence 
of the (subsequence) $f_{(c_j,K_j)}$ to $f$ on $\partial \triangle$, for $j$ large, the ellipses $E_j(c_j)$ lie in any apriori prescribed neighborhood of $E$. Fix such a neighborhood $N$ of $E$ with $N\cap L$ lying to the left of $p$. For
large $j$, we can slide $E_j(c_j)$ to the right in $K_j$ and expand $E_j(c_j)$ to get an ellipse $E_j'$ with a larger scale factor $\rho$, contradicting the
extremality of (the unique) ellipse $E_j(c_j)$ for $K_j,c_j$. 
\medskip

\noindent {\sl Case II: $p \leq T_1(c)$ or $p \geq T_2(c)$}: We may assume $p \leq T_1(c)$. Using the same notation as in Case I; i.e., $E_j(c_j)=f_{(c_j,K_j)}(\partial \triangle)$, we show that 
$\lim_{j\to \infty} a_0(c_j,K_j)=:a$ exists; moreover, in this case, we show that the ellipse 
$$\{a+ \rho(c,K)\{(1,c)\zeta^{-1} +(1,\bar c)\zeta\}: |\zeta|=1\}$$
has top $t=T_1(c)$. To see this, note first that we cannot have a subsequence of the $l_j'$s converging to $\tilde p < T_1(c)$, for by Corollary 9.1 we would then 
get an extremal ellipse for $K,c$ centered to the left of $a_1(c)$. Suppose for the sake of obtaining a contradiction that a subsequence of the $l_j'$s converges
to $\tilde p > T_1(c)$. Form the ellipse $E$ with top at $\tilde p$ and bottom $\tilde b$ on $M$ using the scale factor $\rho(c,K)$; i.e., 
$$E=\{f(\zeta):=\tilde a+  \rho(c,K)\{(1,c)\zeta^{-1} +(1,\bar c)\zeta\}: |\zeta|=1\}$$
where $\tilde a$ is the midpoint of $\tilde p$ and $\tilde b$. Corollary 9.1 implies that $E$ is an extremal ellipse for $K,c$; moreover, by the uniform convergence 
of the (subsequence) $f_{(c_j,K_j)}$ to $f$ on $\partial \triangle$, for $j$ large, the ellipses $E_j(c_j)$ lie in any apriori prescribed neighborhood of $E$. Fix such a neighborhood $N$ of $E$ with $N\cap L$ lying to the right of $T_1(c)$. For
large $j$, we can slide $E_j(c_j)$ to the left in $K_j$ and expand $E_j(c_j)$ to get an ellipse $E_j'$ with a larger scale factor $\rho$, again contradicting the
extremality of $E_j(c_j)$ for $K_j,c_j$. 

\medskip

In case the $x-$axis is a horizontal tangent to $\partial K$ at an endpoint $(0,0)$ or $(2p,0)$ of $L$, we slightly modify the homotopy $H$ so that it is no 
longer a fixed-endpoint homotopy; but we must insure that the approximating sets $K_j$ are convex. Suppose $\partial K$ has a horizontal tangent at $(2p,0)$. Then
we extend this side $L$ to $b:=(2p+2\epsilon,0)$ and connect $b$ to $\partial K$ with a segment emanating from $b$ so that the resulting set $\tilde K \supset K$ is
convex. We now construct a linear midpoint homotopy $H(\alpha,t)$ with $p$ replaced by $p+\epsilon$ on the segment from $(0,0)$ to $b$. In the ``side-by-side''
procedure described using the linear midpoint homotopies $H_j$ for $L_j$, i.e., the choice of $\{\alpha^1_j\}_{j=1,...}$ to form $K_1$, $\{\alpha^2_j\}_{j=1,...}$
to form $K_2$, etc., we now utilize the modified homotopies wherever a horizontal tangent occurs and choose sequences $\{\alpha^k_j\}$ as well as $\epsilon_k
\downarrow 0$ so that the constructed sets $K_j$ belong to ${\mathcal Q}$ with $K_j\downarrow K$. That is, we simply modify the procedure outlined in our previous
linear midpoint homotopy construction and then all the arguments follow as before.

\end{proof}

\noindent {\bf Remark}. It is possible to have a continuous map
$$F_K(c/\zeta)=a(c)+ \rho(c,K)\{(1,c_2,...,c_n)\zeta^{-1} +(1,\bar c_2,...,\bar c_n)\zeta\}$$
to $\CC\PP^n\setminus K$; i.e., the center function $c\to a(c)$ is continuous, but for which the leaves $C,\tilde C$ corresponding to distinct $c,\tilde c\in H_{\infty}$ are not necessarily disjoint in $\CC\PP^n \setminus K$. To see this, let $K=[-1,1]\times [-1,1]$ be a square in $\RR^2$; take $c=[0:1:0]$ and $\tilde c =[0:1:ib]$ for $b>0$ sufficiently small. We choose $\alpha \in [-1,1]$ and define 
$C=f_c(\CC\PP^1)$ where $$f_c(\zeta)=(\frac{1}{2}(\zeta +1/\zeta),\alpha),$$
and we choose $\beta \in [-1+b,1-b]$ and define
$\tilde C=f_{\tilde c}(\CC\PP^1)$ where
$$f_{\tilde c}(\zeta)=(\frac{1}{2}(\zeta +1/\zeta),\frac{ib}{2}(\zeta -1/\zeta)+\beta).$$
If $|\alpha - \beta|>b$, it follows that $C$ meets $\tilde C$ outside of $K$. Note that $\tilde c \to c$ as $b\to 0$. Fix $\alpha$ and define $\beta =\beta(b)$ so that $|\alpha - \beta(b)|>b$ for $b>0$ sufficiently small. Clearly we can define a continuous function 
$$a:H_{\infty}\to (\{0\}\times [-1,1])\bigcup ([-1,1]\cup \{0\})$$ 
with $a(c)=(0,\alpha)$ and $a(\tilde c)=(0, \beta)$.

\vskip5mm

\section{\bf $F_K$ is a homeomorphism, case III.}
\label{sec:FhomeoIII}

\vskip3mm 

We have seen that for a convex body $K$ in $\RR^n$, the additional hypothesis of property ${\mathcal Q}$ (uniqueness of extremal curves) is essentially equivalent to the existence of a unique, well-defined center function $a_0(c)$.  The Lempert theory approximation by $\sigma-$invariant strictly lineally convex domains $D_j \downarrow K$ described in section 5 and used in \cite{blm} guarantees the continuity of the center function. Where one does not have unique extremals, one might still be able to make a canonical choice of center function (e.g., symmetric bodies with respect to 0, where the canonical choice is $a_0(c)\equiv 0$; cf., the remark in section 7).  A natural generalization for a canonical center function would be the barycenter of allowable centers of extremal ellipses.  Precisely, for each $c\in H_{\infty}$, let 
$$A(c):=\{a\in K: a+\rho(c,K)\{(1,c_2,...,c_n)\zeta^{-1} +(1,\bar c_2,...,\bar c_n)\zeta\}$$
$$ \ \hbox{is extremal for} \  K\}$$
and define the center function
$c\mapsto b(c)$, where $b(c)$ is the barycenter of the convex set $A(c)$ in $\RR^n$.

We will show that this gives a foliation in the $\RR^2$ case.  One has to deal with a couple of issues.  First, in prescribing a center function a priori, it is not known that such a center function can be obtained via approximation by strictly linearly convex domains.  To get around this we will show the foliation property directly, using properties of the extremal curves for convex bodies.  We take section 8 as our starting point in this investigation.  Also, property ${\mathcal Q}$ was used in several parts of section 9 to show that the extremal curves were disjoint in $\CC\PP^n\setminus K$.  We will use similar arguments to show that the barycenter construction gives disjoint leaves.

In \cite{blm} it was shown that there exists through each point of $\CC^n\setminus K$ an extremal curve on which $V_K$ is harmonic.  In choosing a specific center function $a_0(c)$ we are essentially throwing away extra leaves in the hope that the remaining leaves give us a foliation.  It turns out that requiring these leaves to be disjoint automatically ensures that they fill out all of $\CC\PP^n\setminus K$.  Using the map $F_K:\CC\PP^n\setminus K_{\rho}\to \CC\PP^n\setminus K$ of (\ref{eqn:fkc}) with a specific center function $a(c)$, the precise statement reads as follows.

\begin{proposition} \label{onto} Suppose $$F(c/\zeta):=a(c)+\rho(c,K)\{(1,c_2,...,c_n)\zeta^{-1} +(1,\bar c_2,...,\bar c_n)\zeta\}$$
is continuous and 1-1; in particular, suppose $c\to a(c)$ is continuous. Then $F$ maps $\CC\PP^n\setminus K_{\rho}$ onto $\CC\PP^n\setminus K$.
\end{proposition}

\begin{proof}  From the construction of the leaves it is clear that they cover every point of $H_{\infty}\subset \CC\PP^n\setminus K$; i.e., $F(H_{\infty})=H_{\infty}$. Thus it suffices to show that $F(\CC^n\setminus K_{\rho})=\CC^n\setminus K$.  Since $\CC^n\setminus K$ is connected, we can do this by verifying that $F(\CC^n\setminus K_{\rho})$ is both open and closed in $\CC^n\setminus K$.  That $F(\CC^n\setminus K_{\rho})$ is open can be seen from the principle of invariance of domain. To show that $F(\CC^n\setminus K_{\rho})$ is closed, let $q\in\overline{F(\CC^n\setminus K_{\rho})}$ (here the closure is taken relative to $\CC^n$) and let $\{q_j\}\subset F(\CC^n\setminus K_{\rho})$ be a sequence converging to $q$.  Now a normal families argument as in  \cite{blm} shows that a subsequence of the leaves through each $\{q_j\}$ converges to a limiting leaf $C$.  This limiting curve is given by $\zeta\mapsto a'+\rho(c')\{(1,c'_2,...,c'_n)\zeta^{-1} +(1,\bar c'_2,...,\bar c'_n)\zeta\}$ for some $c'\in H_{\infty}$ and some $a'$; in particular, $c'=\lim_{j\to \infty} c^{(j)}$ for some $\{c^{(j)}\}$.  But then by continuity of the center function we must have $a'=a(c')$, which in turn implies that 
$$q=\lim_{j\to \infty} q_j = \lim_{j \to \infty} F(c^{(j)}/\zeta_j )=F(c'/\zeta' )$$
for some $\zeta'=\lim_{j\to \infty} \zeta_j$ (since $V_K(q_j)\to V_K(q)$) so that 
$q\in C\subset F(\CC^n\setminus K_{\rho})$.  Thus $F(\CC^n\setminus K_{\rho})$ is closed.

\end{proof}

We now specialize to convex bodies in $\RR^2$. We again use the notation $c=[0:1:c]\in H_{\infty}$. 

\begin{proposition} \label{dfaces} Two extremal curves $C$ and $C'$ are disjoint in $\CC\PP^2\setminus K$ if at least one of the following is true.
\item{(i)} At least one of the extremals is unique $($for its value of $c\in H_{\infty})$; or, if not,
\item{(ii)} If $C$ intersects the parallel edges $X_1$ and $X_2$ then  $C'$ does not intersect both of $X_1$ and $X_2$.
\end{proposition}

\begin{proof} We consider these two cases.
\item{(i)} This follows from the remark at the end of section 8.

\item{(ii)}  Suppose $C\cap X_1=\{z_0\}$ and $C'\cap X_1=\emptyset$.  Form the convex hull $\mathcal{H}$ of $(C\cup C')\cap K$.  We claim that $C$ is a unique extremal for the set $\mathcal{H}$.  Clearly it is an extremal for $\mathcal{H}$ since it is extremal for $K$.  Note that $(C\cup C')\cap X_1 = C\cap X_1=\{z_0\}$.  Now all extremal curves associated to $c$ for $K$ must also intersect $X_1$ and $X_2$, and $C$ is the only one that intersects $X_1$ at $z_0$.  Hence $C$ must be unique for $\mathcal{H}$.  We may now apply case (i) to $C$, $C'$ and $\mathcal{H}$.

\end{proof}

Suppose $X_1$ and $X_2$ are a pair of parallel edges of a convex body $K\subset\RR^2$.  The centers of all  extremal ellipses that touch both edges lie on the parallel line that is equidistant to the lines containing $X_1$ and $X_2$.  By making a suitable coordinate change we may take this line to be the $x$-axis in $\RR^2$. For values of $c\in H_{\infty}$ that give these extremal ellipses, the corresponding allowable centers are intervals on the $x$-axis, i.e., 
$$
A(c)=[\alpha_1(c),\alpha_2(c)]\times\{0\}.
$$
Let $a_1(c)=(\alpha_1(c),0)$ and $a_2(c)=(\alpha_2(c),0)$.

\begin{lemma} \label{ends} Suppose that distinct points $c,c'\in H_{\infty}$ generate possibly nonunique extremal curves  that touch $X_1$ and $X_2$.  Then the particular extremals $$\zeta\mapsto a_1(c)+\rho(c,K)\{(1,c)\zeta^{-1} +(1,\bar c)\zeta\}$$
and
$$\zeta\mapsto a_1(c')+\rho(c',K)\{(1,c')\zeta^{-1} +(1,\bar c')\zeta\}$$ are disjoint in $\CC\PP^2\setminus K$.  The result is also true if $a_1(c),a_1(c')$ is replaced with $a_2(c),a_2(c')$.
\end{lemma}

\begin{proof}  We will prove only the case where we have two nondegenerate leaves, as the other cases follow identical reasoning.   Let $C_{\RR}$ and $C'_{\RR}$ be the inscribed ellipses in $K$ corresponding to the given leaves.  Without loss of generality, we may assume that $\alpha_1(c)\leq\alpha_1(c')$.  Since two ellipses can intersect in 0, 2 or 4 points, we need to show that $C_{\RR}$ and $C'_{\RR}$ intersect in at least three points (counting multiplicities).  If not, then in $\RR^2$ we have at most two points of intersection. But then one can slide the ellipse $C'_{\RR}$ farther to the left to get further extremals for $K,c'$ with centers to the left of $a_1(c')$. This contradicts the definition of $\alpha_1(c')$.

\end{proof}

\begin{theorem} For each $c\in H_{\infty}$, let $b(c)$ be the midpoint of the line segment $A(c)$.  Then 
$$
F(c/\zeta):=b(c)+\rho(c,K)\{(1,c)\zeta^{-1} +(1,\bar c)\zeta\}
$$
is a homeomorphism.

\end{theorem}

\begin{proof} The map $c\to b(c)$ is continuous from Lemma \ref{topcont}. By Proposition \ref{onto} it suffices to show that the extremals $\zeta\mapsto b(c)+\rho(c,K)\{(1,c)\zeta^{-1} +(1,\bar c)\zeta\}$, $c\in H_{\infty}$, are disjoint in $\CC\PP^2\setminus K$.  By Proposition \ref{dfaces}, we need only consider pairs $c,c'\in H_{\infty}$ whose extremals touch the same pair of parallel edges, and show that the extremals centered at $b(c)$ and $b(c')$ are disjoint in $\CC\PP^2\setminus K$.

Fixing $c$ and $c'$, without loss of generality we may assume as in Lemma \ref{ends} that
$$
A(c)=[\alpha_1(c),\alpha_2(c)]\times\{0\}, \hbox{ and } A(c')=[\alpha_1(c'),\alpha_2(c')]\times\{0\}.
$$

We show that the real extremal ellipses centered at $b(c)$ and $b(c')$, where
$$
b(\cdot)=(\beta(\cdot),0):=\bigl(\frac{1}{2}(\alpha_1(\cdot)+\alpha_2(\cdot)),0\bigr)
$$
have more than two points of intersection.  Let $j$ be the number of points of intersection of these ellipses.  Without loss of generality, we assume that $\alpha_2(c)-\alpha_1(c)\geq\alpha_2(c')-\alpha_1(c')$ and that $\beta(c)\geq\beta(c')$.  Then we see that the extremals centered at $a_2(c)$ and $a_2(c')$ lie on the same line but are farther apart than those centered at $b(c)$ and $b(c')$.  Thus they can have at most $j$ points of intersection.  But if $j\leq 2$ this contradicts Lemma \ref{ends}.

\end{proof}

\noindent {\bf Remark}. In higher dimensions, the collection of centers, $A(c)$, may not vary continuously with $c$ (e.g. in the Hausdorff metric). For example, let
$K\subset\RR^3$ be the convex hull of the set
$$
\{(x_1,x_2,x_3):x_1\in [-1,1],\ x_2=x_3=0\}\cup \{(x_1,x_2,x_3):x_3=1, \ x_1^2+x_2^2\leq 1\}.
$$
All extremal ellipses centered at
$(0,0,1)$ with major axis not parallel to the $x_1$ axis are unique; i.e., for such $c$, $A(c)=\{(0,0,1)\}$. But if an extremal ellipse centered at $(0,0,1)$ other than the unit circle has major axis parallel to the $x_1$-axis, then it can slide along the $x_3$-axis, so that, for such $c$, $A(c)$ is a nondegenerate segment with one endpoint at $(0,0,1)$. Thus $c\mapsto A(c)$, and hence $c\mapsto b(c)$, is not continuous.
\medskip

We conclude this section with a simple example of a foliation associated to a non-convex body. We will say that a compact set $K$ has the {\it foliation property} if there is a foliation of $\CC^n\setminus K$ by one-dimensional complex varieties such that $V_K$ is harmonic on each leaf of the foliation. We make the following observation:
{\sl Let $K\subset\CC^n$ be a compact set, and suppose there exists a polynomial map $P:\CC^n\to\CC^n$ such that $P(K)$ has the foliation property, $P^{-1}\circ P(K)=K$, and 
\begin{equation}
\label{eqn:proper}
V_K(z)=\frac{1}{\deg P}V_{P(K)}(P(z)).
\end{equation} 
Then $K$ has the foliation property, with leaves given by $P^{-1}(L)$ for each leaf $L$ of the foliation associated to $P(K)$.}

\medskip

\noindent{\bf Example.}  Let $z_1,z_2$ be the coordinates in $\CC^2$ with $x=\Re z_1$ and $y=\Re z_2$ the corresponding coordinates in $\RR^2$.  Let $K\subset\RR^2$ be the annular region bounded by the circle $x^2+y^2=4$ and the ellipse $x^2+4y^2=1$.  Consider the map defined by $P(z_1,z_2)=(z_1^2,z_2^2)$. This is a proper map from $\CC^2$ onto $\CC^2$ and thus by Theorem 5.3.1 of \cite{klimek}, (\ref{eqn:proper}) holds. Then $P(K)$ is the convex quadrilateral bounded by the coordinate axes together with the lines $x+y=4$ and $x+4y=1$.  By Theorem 7.1, $P(K)$ is a convex set in $\RR^2$ with the foliation property, whose $P$-preimage is $K$.  Thus by the above observation, $K$ has the foliation property.

We take a closer look at the leaves for $K$; i.e., the pullbacks under $P$ of leaves for $P(K)$.  For instance, from section 6, the complex line $L:z_2=z_1+\frac{1}{4}$ is a (degenerate) leaf of the foliation for $P(K)$.  This pulls back to the complex hyperbola $z_2^2-z_1^2=\frac{1}{4}$. Similarly, a degenerate leaf $L$ for $P(K)$ given by a complex line of the form $z_2=-az_1+b$ with $1/4 \leq b \leq 2$ and $a=b/4$ pulls back to an ellipse $P^{-1}(L)$ such that $E:=\RR^2\cap P^{-1}(L)$ is a real ellipse.  Note that in contrast to the convex case, $E$ is not necessarily contained in $K$. Nondegenerate ellipses $L$ pull back to varieties given by degree 4 polynomials.

We may approximate $P(K)$ from above by a decreasing sequence of relatively compact, strictly lineally convex domains $\{D_j\}$ in $\CC^n$ that are invariant under conjugation, $\sigma(D_j)=D_j$, so that $P(K)=\cap_j D_j$. The pre-images $\{P^{-1}(D_j)\}$ each have the foliation property, are conjugation invariant, and decrease to $K$: $\cap_j P^{-1}(D_j)=K$. By the remark at the end of section 5, for $j$ large the sets $P^{-1}(D_j)$ cannot be weakly lineally convex. Thus we have exhibited examples of non-convex bodies in $\RR^n$ as well as non-lineally convex domains in $\CC^n$ which have the foliation property.

 \vskip 5mm

\section{\bf Robin indicatrix for convex bodies and symmetrization.}
\label{sec:ris}

\vskip3mm

In this section, we give a relationship between the Robin indicatrix of an arbitrary convex body in $\RR^n$ and its so-called symmetrization. Let $K$ be a non-pluripolar compact set in $\CC^n$; let $V_K$ be the Siciak-Zaharjuta extremal function; let 
$$\rho_K(w):=\limsup_{|\lambda|\to \infty, \ \lambda \in \CC}[V_K^*(\lambda w)- \log {|\lambda|}]$$
be the logarithmically homogeneous Robin function for $K$; and recall from (\ref{krho}) that $K_{\rho}=\{w\in \CC^n: \rho_K(w)\leq 0\}$. This corresponds to the complement of the set $R(D)$ if $K=\bar D$ is the closure of a strictly lineally convex domain. Then the lines 
$w=v/\zeta $ for $v\in \partial K_{\rho}$ and $\zeta \in \triangle$ are lines on which 
\begin{equation}
\label{eqn:rob1}\rho_K(w)=\rho_K(v/\zeta )=\rho_K(v)-\log {|\zeta|}=-\log {|\zeta|}
\end{equation} 
since $\rho_K(v)=0$. For $K$ a {\it symmetric} convex body in $\RR^n$, the Robin exponential map $F_K:\CC\PP^n\setminus K_{\rho} \to \CC\PP^n \setminus K$ as described in section 5 (see (\ref{eqn:fk})) restricted to $\CC^n\setminus K_{\rho}$ is given by 
\begin{equation}
\label{eqn:baran}
z=F_K(w)=F_K(v/\zeta)=v/\zeta+\bar v\zeta.
\end{equation}
Note that, up to a factor of $2$, this is the same as Baran's parameterization in \cite{bar}. Explicitly, we have from (\ref{eqn:rob1}) that if $\rho_K(v)=0$
$$\zeta =e^{-\rho_K(w)}e^{i\theta(w)} \ \hbox{and} \ v=w e^{-\rho_K(w)}e^{i\theta(w)}$$
so that
$$z=F_K(w)=w+ \frac{\bar w}{e^{2\rho_K(w)}}=e^{-\rho_K(w)}\bigl [e^{\rho_K(w)}w+e^{-\rho_K(w)}\bar w\bigr ]. $$
This is a ``real'' map; i.e., $w\in \RR^n\setminus K_{\rho}$ if and only if $z=F_K(w)\in \RR^n\setminus K$: for let $\rho_K(w)=\log {t}>0$. Then $z=\frac{1}{t}[tw +\bar w/t]$ so that $z\in \RR^n$ if and only if $tw +\bar w/t=t\bar w +w/t$; i.e., if and only if $t(w-\bar w)=(w-\bar w)/t$; i.e., if and only if $w=\bar w$.

Note that if we consider points with $|\zeta|=1$, then the image of $\partial K_{\rho}$, the zero-level set of $\rho_K$, is $K$, the zero-level set of $V_K$. Precisely, we have 
$$F_K(\partial K_{\rho})=\{2 \, \Re \, w: w\in \partial K_{\rho}\}=K$$
since $\rho_K(w)=0$ for such points. That is, modulo scaling by a factor of $2$, the projection of $\partial K_{\rho}$ onto $\RR^n$ recovers our set $K$ (equivalently, the intersection of $\partial K_{\rho}$ with $\RR^n$ is $\frac{1}{2}K$). For example, the interval $K=[1,1] \subset \RR \subset \CC$ has Robin function $\rho_K(w)=\log {|w|} + \log 2$ so that $K_{\rho}=\{w: |w|\leq 1/2\}$ and $\partial K_{\rho}$ projects to $[-1/2,1/2]$.

\begin{corollary}
	 For a symmetric convex body $K$ in $\RR^n$,  
for each $t>1$,
$$K(t):=\{x\in \RR^n: V_K(x) \leq \log t\}$$
is a homothety of $K$.
\end{corollary}

\begin{proof} Under the map $F_K$, level sets of $\rho_K$ map to level sets of $V_K$ (by construction). Thus 
$$\partial K(t) \ \hbox{(in $\RR^n$)} \ = \{x\in \RR^n: V_K(x) = \log t\}=F_K(\{w: e^{\rho_K(w)}=t\})\cap \RR^n$$
$$= \{z=\frac{1}{t}[tw +\bar w/t]: e^{\rho_K(w)}=t\}\cap \RR^n=\frac{1}{t} \{(t+1/t)w: w\in \RR^n, \ e^{\rho_K(w)}=t\}$$
$$=\frac{1}{t} \{(t^2+1)(w/t): w\in \RR^n, \ \rho_K(w/t)=0\}=(t+1/t)\bigl[\partial K_{\rho}\cap \RR^n\bigr]$$
$$=\frac{1}{2}(t+1/t)\partial K.$$
\end{proof}

Now suppose $K$ is an arbitrary convex body in $\RR^n$ (not necessarily symmetric). The mapping $F_K$ is still well-defined as a map from $\CC^n\setminus K_{\rho}$ to $\CC^n \setminus K$:
$$z=F_K(w)=F_K(v/\zeta)=a_0(v)+ v/\zeta +\bar v \zeta,$$
but since we don't know an explicit formula for the center function $a_0(v)$, it may be difficult to write down an explicit formula for $F_K(w)$ if $v\in \partial K_{\rho}$ in terms of $w$. However, we can do the following: form $K_{\rho}$ and define the ``forgetful center map'' 
$$z=\tilde F_K(w)=\tilde F_K(v/\zeta)=v/\zeta+\bar v \zeta$$
on $\CC^n\setminus K_{\rho}$ and again considering points with $|\zeta|=1$ we let $\tilde K:=\tilde F_K(\partial K_{\rho})$. Since (\ref{eqn:rob1}) holds, as before, if $v\in \partial K_{\rho}$,
$$z=\tilde F_K(w)= e^{-\rho_K(w)}\bigl [e^{\rho_K(w)}w+e^{-\rho_K(w)}\bar w\bigr ] $$
and 
$$\tilde K=\{2\Re w: w\in \partial K_{\rho}\}$$
since $\rho_K(w)=0$ for such points.

\vskip6pt
\noindent (1) {\sl $\tilde K$ is a symmetric convex body in $\RR^n$}.
\vskip4pt
The symmetry is obvious (e.g., since $K_{\rho}$ is balanced). To show the convexity, first of all, from \cite{momm}, the sublevel sets of $V_K$ intersected with $\RR^n$ are convex; i.e., for each $t>1$,
$$K(t)=\{x\in \RR^n: V_K(x) \leq \log t\}$$
is convex. Now using the {\it definition} of $\rho_K$, we show that the sublevel sets of $\rho_K$ intersected with $\RR^n$ are convex: given $t>1$, if $\rho_K(x), \ \rho_K(y) \leq  \log t$, then given $\epsilon >0$ we have 
$$\max [V_K(\lambda x) - \log {\lambda} -\log t,V_K(\lambda y) - \log {\lambda} -\log t] < \epsilon$$
for $\lambda\in \RR^n$ with $\lambda > M(\epsilon)$. Then $\lambda x, \lambda y \in K(\lambda t e^{\epsilon})$. By convexity, 
$a\lambda x  +(1-a)\lambda y \in K(\lambda t e^{\epsilon})$ for $0\leq a \leq 1$; i.e., 
$$V_K(a\lambda x  +(1-a)\lambda y) -\log {\lambda} =V_K(\lambda (ax+(1-a)y))-\log {\lambda}$$
$$ < \log t +\epsilon$$
for $\lambda > M(\epsilon)$. This says that $\rho_K(ax +(1-a)y) \leq \log t$. It follows that $\tilde K$ is convex.

\vskip6pt

\noindent (2) {\sl $V_{\tilde K} (v/\zeta  +\bar v \zeta)=-\log {|\zeta|}$
for $v\in \partial K_{\rho}$ and $\zeta\in \triangle\setminus \{0\}$; that is, $\tilde K_{\rho}=K_{\rho}$ and the leaves in the foliation of $\CC^n\setminus \tilde K$ are the origin-centered translates of the leaves in the foliation of $\CC^n\setminus K$.}
\vskip4pt

The fact that $\tilde K_{\rho}=K_{\rho}$ follows since $F_K$ and $\tilde F_K$ differ by the (bounded) center function $a_0(v)$; the rest follows from Baran \cite{bar} (see (\ref{eqn:baran})).

\vskip6pt
\noindent (3) {\sl $\tilde K$ is the ``symmetrization'' $K_{sym}$ of $K$}. 
\vskip4pt

For a convex set $E$ in $\RR^n$, the symmetrization $E_{sym}$ is the difference set
$$E_{sym}:=\frac{1}{2}(E - E) := \{\frac{1}{2}(x-y):x,y\in E\}.$$
Note that if $E$ is symmetric, then $E_{sym}=E$. It follows easily from the definition of $E_{sym}$ that 
\vskip3pt
(a) if $E$ is an interval, $E_{sym}$ is a symmetric interval in the same direction as $E$ and with the same length;

(b) if $E$ is an ellipse, $E_{sym}$ is a translated, symmetric ellipse with the same eccentricity and orientation as $E$.
\vskip3pt

Since $F_K^{-1}(a_0(v)+ v/\zeta  +\bar v \zeta)= v/\zeta$ and 
$\tilde F_K(v/\zeta)=v/\zeta  +\bar v \zeta$, 
we have 
$$( \tilde F_K \circ F_K^{-1})(a_0(v)+ v/\zeta  +\bar v \zeta)=v/\zeta  +\bar v \zeta.$$
Thus for $|\zeta|=1$ this maps $K$ onto $\tilde K$ by translating each ellipse $e(v):=\{a_0(v)+ v/\zeta  +\bar v \zeta:|\zeta|=1\}$ so that its center is at the origin; i.e., $$(\tilde F_K \circ F_K^{-1})(e(v))= \{v/\zeta  +\bar v \zeta:|\zeta|=1\}.$$ From (2), (a) and (b) it follows that $K_{sym}=\tilde K$.

In contrast to Corollary 3.2, we prove the following. 

\begin{proposition} Let $K\subset \RR^n$ be a convex body. The set  $K_{\rho}$ is not strictly convex $($in $\CC^n)$. 
\end{proposition}

\begin{proof} First of all, if $K$ is the real ball 
$$K:=\{ x=(x_1,...,x_n)\in \RR^n: |x_1|^2+\cdots + |x_n|^2 \leq 1\},$$
setting $z_j=x_j+iy_j$ it is known that 
$$V_K(z_1,...,z_n)= \frac{1}{2} \log h(|z_1|^2+ \cdots + |z_n|^2 +|z_1^2 + \cdots + z_n^2-1|)$$
where $h(\zeta) =\zeta + (\zeta^2-1)^{1/2}$ with the square root chosen so that $h(t) >1$ for $t>1$ (cf., \cite{bar}). Then 
$$\rho_K(z_1,...,z_n)=  \frac{1}{2} \log [2(|z_1|^2+ \cdots + |z_n|^2 +|z_1^2 + \cdots + z_n^2|)]$$
and we thus have
$$K_{\rho}=\{z=(z_1,...,z_n): |z_1|^2+ \cdots + |z_n|^2 +|z_1^2 + \cdots + z_n^2| \leq 1/2\}.$$
It is easy to see that at any point $p\in \partial K_{\rho}\cap \RR^n$ there is a nondegenerate real line segment $S_p$ contained in $\partial K_{\rho}\cap T_p(\partial K_{\rho})$ (e.g., at $p=(1/2,0,...,0)$, all points of the form $(1/2,ti/2,0,...,0)\in \partial K_{\rho}\cap T_p^{\CC}(\partial K_{\rho})$ for $t\in [-1,1]$). Thus $K_{\rho}$ is not strictly convex at any point $p\in \partial K_{\rho}\cap \RR^n$. 

An immediate corollary of this calculation is that the same result holds for a centrally symmetric ellipsoid
$$E:=\{ x=(x_1,...,x_n)\in \RR^n: |a_1x_1|^2+\cdots + |a_nx_n|^2 \leq 1\};$$
i.e., at any point $p\in \partial E_{\rho}\cap \RR^n$ there is a nondegenerate real line segment $S_p$ contained in $\partial E_{\rho}\cap T_p(\partial E_{\rho})$. To see this, consider the linear map $z\to P(z):=(a_1z_1,...,a_nz_n)$ which maps $E$ onto $K$. Using (\ref{eqn:proper}), we have 
$$V_E(z)=V_K(P(z))= V_K(a_1z_1,...,a_nz_n)$$
so that $\rho_E(z) = \rho_K(a_1z_1,...,a_nz_n)$ as well.

Next, take any strictly convex, symmetric convex body $K\subset \RR^n$ and fix $\tilde p\in \partial K$. Consider an inscribed ellipsoid $E_i$ to $K$ at $\tilde p$. Then $(E_i)_{\rho}\subset K_{\rho}$ and the point $p=\tilde p/2 \in \RR^n$ is a common boundary point of these sets. Denote the (common) tangent space at $p$ by $H$. Applying the previous paragraph's result to $(E_i)_{\rho}$, we have that $H\cap \partial (E_i)_{\rho}$ contains a nondegenerate real line segment $S_p$. Suppose $K_{\rho}$ is strictly convex at $p$. In particular, $H\cap \partial K_{\rho}=\{p\}$. But then
$$S_p\subset H\cap \partial (E_i)_{\rho}\subset H\cap \partial K_{\rho} =\{p\},$$
contradicting the fact that $S_p$ is nondegenerate.

Now suppose $K$ is a symmetric convex body in $\RR^n$ which is not strictly convex. We saw that the intersection of $\partial K_{\rho}$ with $\RR^n$ is $\frac{1}{2}K$; hence $K_{\rho}$ is not strictly convex. 

Finally, if $K$ is any convex body in $\RR^n$, form the symmetrization $K_{sym}$ of $K$. Then $K_{\rho}=(K_{sym})_{\rho}$ and the result follows from the symmetric case.

\end{proof}

Note that the Robin indicatrix does not uniquely determine the convex body $K$. In a future work we hope to describe invariants associated to $K$ and its indicatrix in the spirit of \cite{lem3} which will uniquely determine this set.

 \vskip 5mm


\bigskip



\end{document}